\PassOptionsToPackage{unicode}{hyperref}
\PassOptionsToPackage{naturalnames}{hyperref}

\documentclass[a4paper,11pt]{article}

\usepackage{etex}
\usepackage{graphicx}
\usepackage{amsmath}
\usepackage{amssymb}
\usepackage{amsthm}
\usepackage{empheq}
\usepackage{multirow}
\usepackage{makecell}
\usepackage{pict2e}
\usepackage{booktabs}
\usepackage[figuresright]{rotating}
\usepackage{bbm}

\usepackage{algorithm}      
\usepackage{pgfplots}
\pgfplotsset{compat=1.17}
\usepackage{listings}		
\usepackage[noend]{algpseudocode}
\usepackage{changepage}     
\usepackage{color}              
\definecolor{myred}{gray}{0} 
\definecolor{light-gray}{gray}{0.8}

\usepackage[caption=false]{subfig}
\usepackage{epstopdf}

\usepackage{graphics} 
\usepackage{epsfig} 
\usepackage[ plainpages = false, pdfpagelabels,
	     bookmarks,
	     bookmarksopen = true,
	     bookmarksnumbered = true,
	     breaklinks = true,
	     linktocpage,
	     pagebackref,         
	     colorlinks = false,  
         hypertexnames=false,   
	     linkcolor = blue,
	     urlcolor  = blue,
	     citecolor = red,
	     anchorcolor = green,
	     hyperindex = true,
	     hyperfigures
	     ]{hyperref}
\usepackage{pdflscape}
\usepackage{overpic}
\makeatletter
\newenvironment{subfigures}
 {\begin{minipage}{\columnwidth}\def\@captype{figure}\centering}
 {\end{minipage}}
\makeatother

\newtheorem{theorem}{Theorem}[section]

\newtheorem{assumption}[theorem]{Assumption}

\theoremstyle{remark}
\newtheorem{remark}[theorem]{Remark}

\theoremstyle{definition}

\newcommand{\vect}[1]{\textrm{\boldmath${#1}$}} 
\newcommand{\Vect}[1]{{\bf #1}} 
\newcommand{\vel}{u} 
\newcommand\var{w}   
\newcommand\Var{W}   

\newcommand\Fig{figure }

\let\originalleft\left
\let\originalright\right
\renewcommand{\left}{\mathopen{}\mathclose\bgroup\originalleft}
\renewcommand{\right}{\aftergroup\egroup\originalright}

\DeclareMathOperator*{\argmin}{arg\,min}

\binoppenalty=\maxdimen 
\relpenalty=\maxdimen   

\numberwithin{equation}{section}    
\date{\today}             

\title{Hierarchical Micro-Macro Acceleration for Moment Models of Kinetic Equations}
\author{
Julian Koellermeier\footnote{Corresponding author, email address {\tt j.koellermeier@rug.nl}} \footnote{Bernoulli Institute, University of Groningen} \footnote{Department of Computer Science, KU Leuven},
Hannes Vandecasteele\footnotemark[\value{footnote}]
}

\begin{document}

\maketitle


\begin{abstract}
Fluid dynamical simulations are often performed using cheap macroscopic models like the Euler equations. For rarefied gases under near-equilibrium conditions, however,  macroscopic models are not sufficiently accurate and a simulation using more accurate microscopic models is often expensive. In this paper, we introduce a hierarchical micro-macro acceleration based on moment models that combines the speed of macroscopic models and the accuracy of microscopic models. The hierarchical micro-macro acceleration is based on a flexible four step procedure including a micro step, restriction step, macro step, and matching step. We derive several new micro-macro methods from that and compare to existing methods. In 1D and 2D test cases, the new methods achieve high accuracy and a large speedup. 
\end{abstract} 

\medskip\noindent

\section{Introduction}
\label{sec:1}

The evolution of a rarefied gas is modeled using the mass density function $f(t,x,c)$ which follows the Boltzmann transport equation \cite{Cercignani1994}. We first focus on the one-dimensional case and discuss the extension to multiple spatial and velocity dimensions in section \ref{sec:7_mM2DHME}. The Boltzmann equation in 1D reads
\begin{equation} \label{e:BTE}
    \partial_t f(t,x,c) + c \partial_x f(t,x,c) = S(f),
\end{equation}
where $t \in \mathbb{R}_+$ denotes the time, $x \in \mathbb{R}$ the position in physical space, and $c \in \mathbb{R}$ the microscopic velocity, respectively. The right-hand side collision operator $S(f)$ can be modeled in different ways. In this paper, we consider a simple BGK operator \cite{Bhatnagar1954}
\begin{equation} \label{e:BGK}
    S(f) = - \frac{1}{\epsilon} \left( f - f_M\right),
\end{equation}
modeling relaxation with relaxation time $\epsilon \in \mathbb{R}_+$ towards the local Maxwellian $f_M(t,x,c)$, typically given by
\begin{equation} \label{e:Maxwellian}
    f_M(t,x,c) = \frac{\rho(t,x)}{\sqrt{2 \pi \theta(t,x)}} \exp\left( -\frac{\left(c-u(t,x)\right)^2}{2 \theta(t,x)} \right).
\end{equation}
Other version of the collision operator exist and can be found in \cite{Cercignani2000,Shakhov1968}. Note that large relaxation times starting from $\epsilon \approx 1$ lead to pertaining deviations from the local Maxwellian and so-called non-equilibrium effects.

By integration over velocity space the macroscopic variables density $\rho(t,x)$, velocity $\vel(t,x)$, and temperature $\theta(t,x)$ can be extracted from the distribution function $f(t,x,c)$:
\begin{eqnarray} \label{e:macroscopic_quantities}
    \rho(t,x) & = & \int_{\mathbb{R}} f(t,x,c) \, dc, \\
    \rho(t,x) \vel(t,x) & = & \int_{\mathbb{R}} c f(t,x,c) \, dc, \\
    \rho(t,x) \theta(t,x) & = & \int_{\mathbb{R}} \left| c-\vel \right|^2 f(t,x,c) \, dc. \label{e:macroscopic_quantities_end}
\end{eqnarray}

The collision operator \eqref{e:BGK} is designed such that it conserves mass, momentum, and energy during collisions. This means that conservation laws can be derived by multiplying the Boltzmann equation \eqref{e:BTE} with monomials $\left( 1, c, c^2 \right)$ and integrating over the microscopic velocity space yielding the well-known Euler equations. The equations in primitive variables and non-conservative form read
\begin{equation}\label{decoupled-scheme}
        \partial_t \left( \begin{array}{c}
        \rho \\
        \vel \\
        \theta \\
      \end{array} \right) +
      \left( \begin{array}{ccc}
         \vel & \rho & 0\\
        \frac{\theta}{\rho} & \vel & 1 \\
        0 & 2 \theta & \vel \\
      \end{array} \right)
       \partial_x \left( \begin{array}{c}
        \rho \\
        \vel \\
        \theta \\
      \end{array} \right) +
      \left( \begin{array}{c}
        0 \\
        0 \\
        \partial_x q \\
      \end{array} \right)= \left( \begin{array}{c}
        0 \\
        0 \\
        0 \\
      \end{array} \right),
\end{equation}
where the heat flux $q(t,x) = \int_{\mathbb{R}} f(t,x,c) \left( c - \vel \right)^3 \, dc$ is either neglected or to be modeled by some microscopic equations. The Euler equations without heat flux yield macroscopic conservation laws that are accurate for flow situations close to equilibrium. This correlates with a small Knudsen number, indicated by $\epsilon \ll 1$ in our model \eqref{e:BTE} and \eqref{e:BGK}. The flow is then relaxing to the Maxwellian \eqref{e:Maxwellian} and is described solely by the macroscopic variables $\rho, \vel, \theta$. 

For large deviations from equilibrium, i.e.. when $\epsilon \approx 1$, additional microscopic variables need to be taken into account to accurately model the non-equilibrium that cannot be expressed by a Maxwellian. Several approaches for these so-called micro models exist. Apart from stochastic models based on the movement of computational particles \cite{Degond2011b,Garcia1999}, moment models are one way to derive physically accurate hierarchical PDE systems that use an extended set of variables, the so-called moments \cite{Struchtrup2006,Torrilhon2016}. Several examples for moment models will be introduced in Section \ref{sec:3}. 

Moment models are based on an expansion of the distribution function using basis polynomials and basis coefficients, which are the moments. A projection onto test functions then leads to PDE systems of the following general form 
\begin{equation}
\label{e:moment_system}
    \partial_t \var  + \Vect{A}\left( \var \right) \partial_x \var  = - \frac{1}{\epsilon}\vect{S}(\var),
\end{equation}
where the variable vector $\var \in \mathbb{R}^M$ contains $M \in \mathbb{N}$ moment variables, typically including $\rho, \vel, \theta$, the system matrix $\Vect{A}\left( \var \right)$ depends on those moments and models the transport properties, while the right-hand side relaxation term $-\frac{1}{\epsilon}\vect{S}(\var)$ models the collisions and typically contains a kernel according to the conservation laws from above.

For small relaxation times, equation \eqref{e:moment_system} becomes stiff and explicit time stepping schemes require very small time steps. Implicit schemes are possible \cite{Mieussens2000}, but are difficult to generalize and to extend to higher order. Other possibilities are splitting schemes, where higher order is not trivial to achieve \cite{Tcheremissine2001}, or IMEX schemes \cite{Pareschi2005}. Another approach is the Heterogeneous Multiscale Method (HMM) \cite{E2007}, which makes use of one (potentially implicit) model for the slow modes and a second model for the fast modes. The Projective Integration (PI) method \cite{Gear2003,Melis2017} uses a comparable technique, performing a number of small inner iterations to damp the fast modes and a subsequent large extrapolation step to evolve the slow modes. The PI method can be interpreted as a HMM version \cite{Maclean2015}. Recently, the PI method has been extended to multiple relaxation times and moment models \cite{Melis2019,Koellermeier2021}. The decomposition into a macroscopic model for the slow modes and a microscopic model for the fast modes was also used to speed up computation of stiff problems in stochastic simulations \cite{Debrabrant2017,Samaey2011}. However, the decomposition in a deterministic macro model and a stochastic micro model can cause difficulties when coupling the models, a disadvantage that can be mitigated using moment models.

In this paper, we propose a hierarchical micro-macro acceleration to obtain fast, yet accurate solutions of the full model by adapting the strategy from \cite{Debrabrant2017} for moment models. The hierarchical micro-macro method achieves an acceleration of the moment simulation in four steps: 
\begin{itemize}
    \item[(1)] microscopic simulation of the moment equations~\eqref{e:moment_system} with a small time step,
    \item[(2)] restriction of the microscopic moment vector to the slow, macroscopic moments, 
    \item[(3)] macroscopic simulation of the restricted moments, 
    \item[(4)] matching to reconstruct a new set of microscopic moments that are consistent with the new macroscopic moments.
\end{itemize}

The goal of the hierarchical micro-macro acceleration is to overcome the severe time step constraint of standard explicit schemes to accelerate the simulations while at the same time achieving high accuracy. The use of hierarchical moment models for the micro and macro models will lead to a large degree of flexibility so that many existing models can be compared to the hierarchical micro-macro acceleration. Special attention has to be given to the matching step, where we opt for a standard $L^2$ matching in this paper. As will be shown, this has the benefit of a cheap, explicit solution of the matching step. Extensions to more advanced matching operators are possible, see \cite{Vandecasteele2020}. The different models will then be tested for a standard 1D test case and a 2D application case. The results show that the new micro-macro methods have a good accuracy together with significant speedup of up to a factor of 250 with respect to a full micro solution.

The rest of the paper is structured as follows: In Section \ref{sec:2}, we propose the hierarchical micro-macro acceleration adapted to the use with general moment models. Examples for new moment models and methods are given in Section \ref{sec:3}, while existing models are rewritten in the hierarchical micro-macro acceleration setting in Section \ref{sec:3b}. Details of the implementation and an analysis of the matching step follow in Section \ref{sec:4}. Numerical results for standard test cases are shown in Section \ref{sec:Results}. The paper ends with a conclusion and further work.    
\section{Hierarchical micro-macro acceleration for moment models}
\label{sec:2}
In this section we introduce the hierarchical micro-macro acceleration for moment models based on the work in \cite{Vandecasteele2020}. The hierarchical micro-macro acceleration consists of four steps and several choices can be made for each of the steps such that it allows for the derivation of new methods and can be compared to existing methods as outlined in detail in the next sections.

Let $\var^n \in \mathbb{R}^{M}$ be the set of moments obtained after $n$ steps of the micro-macro method. One full time step of the hierarchical micro-macro acceleration consists of the following four steps:
\begin{itemize}
    \item[(1)] Microscopic step: we simulate a single small time step $\delta t$ with a microscopic model $g$ to obtain the intermediate micro solution $\var^{n,*} \in \mathbb{R}^{M}$:
    \begin{equation*}
    \var^{n,*} = g(\var^{n}, \delta t).
    \end{equation*}
    \item[(2)] Restriction: we select the first $L \leq M$ moments in the intermediate micro solution $w^{n,*}$ to compute the intermediate macro solution $\Var^{n,*} \in \mathbb{R}^{L}$:
    \begin{equation*}
    \Var^{n,*} = \var^{n,*}_{1:L}.
    \end{equation*}
    \item[(3)] Macroscopic step: we simulate a single large time step $\Delta t$ with some macroscopic model $\mathcal{G}$ to compute the macro solution $\Var^{n+1} \in \mathbb{R}^{L}$:
    \begin{equation*}
         \Var^{n+1} = \mathcal{G}\left(\Var^{n,*}, \Delta t\right).
    \end{equation*}
    \item[(4)] Matching: we reconstruct the new micro solution $\var^{n+1} \in \mathbb{R}^{M}$ that matches the lowest $L$ moments with the macro solution $\Var^{n+1}$, but has the smallest distance to the intermediate micro solution $\var^{n,*}$ measured in some pseudo metric $d$.
    
    
\end{itemize}

\begin{figure}[htb!]
    \centering
    \includegraphics[width=0.75\linewidth]{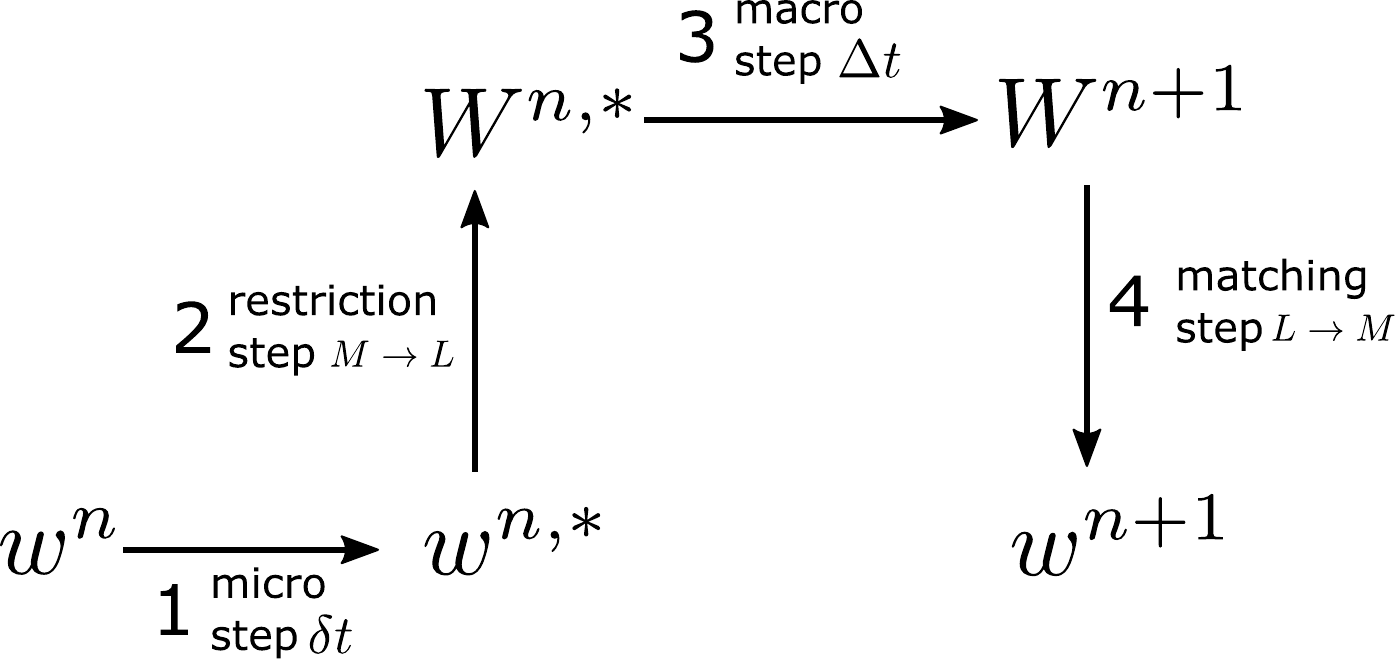}
    \caption{Outline of the hierarchical micro-macro acceleration including the micro step (1), restriction step (2), macro step (3), and matching step (4).}
    \label{fig:mM_sketch}
\end{figure}

The general procedure for one full iteration of the hierarchical micro-macro acceleration proposed in this work is outlined in \Fig \ref{fig:mM_sketch}. Below we will discuss the single steps in more detail.\\


The micro step (1) uses an accurate micro model $g$ containing $M$ microscopic variables, e.g., moments or expansion coefficients of the distribution function. It is typically a PDE model for which standard time stepping schemes are implemented. The time step $\delta t$ is chosen such that it fulfills any stability constraints originating from a possibly stiff right-hand side. This small time step ensures that variables on the micro scale relax quickly towards equilibrium but do not become unstable. 

The restriction step (2) takes the intermediate solution of the micro model and reduces the microscopic variables to a smaller set of macroscopic variables $\Var \in \mathbb{R}^L$ with $L<M$ that evolve on a macroscopic scale, e.g., the conserved quantities mass, momentum, and energy.

The macro step (3) uses a simpler macro model $\mathcal{G}$ to model the dynamics of the macroscopic variables $\Var$. It is typically a PDE model, e.g. the Euler equations, that can be solved by standard time stepping techniques but it can even be a simple extrapolation in time, leading to the Projective Integration method, see \cite{Koellermeier2021}. As the microscopic scales have been eliminated, a time step size $\Delta t$ following a standard CFL condition is sufficient for stability. In practice, this leads to significantly larger time step sizes than for the micro step.

The matching step (4) ensures the consistency of the macro solution with the new micro solution and requires some further explanation. It computes the microscopic variables using the information of the intermediate micro step and the new macro step. The two main principles are: (I) consistent moments of the distribution functions and (II) minimal distance to the prior micro solution, measured in a problem specific norm, for example a standard $L^2$ norm or the Kullback-Leibler divergence. We therefore reconstruct a new set of microscopic variables $\var^{n+1} \in \mathbb{R}^{M}$ that is consistent with the macro variables $\Var^{n+1}$. For this step, we work on the level of mass density functions, and define the macroscopic operator $M_L(f)$ that returns the set of macroscopic moments. Our matching step proceeds by picking the mass density function $f^{n+1}$ such that,
\begin{equation}\label{eq:general_matching}
    f^{n+1} = \underset{f   \in  V(\Var^{n+1})}{\argmin} \ d(f, f^{n,*}),
\end{equation}
for some (pseudo) metric $d$. Here, $V(\Var^{n+1})$ is the function space containing all density functions where the first $L$ moments are constrained by the macroscopic variables $\Var^{n+1}$. Additionally, $f^{n,*}$ is the density function corresponding to the microscopic moments $\var^{n,*}$ after the microscopic step. Finally, $\var^{n+1}$ is defined as the microscopic moments obtained from $f^{n+1}$.

In the following sections, we give examples of three new and two existing methods that can be derived from and compared to the hierarchical micro-macro acceleration above. They differ in the execution of steps 1-4 above. As the restriction and matching step both depend on the choice of the macro model, we will first discuss the micro and macro model, before detailing the restriction and matching steps.

\subsection{Consistency of micro-macro acceleration}
\label{sec:Consistency_proof}
For new methods based on the hierarchical micro-macro acceleration scheme to make sense, we need to show general consistency of this framework. Consistency means that the solution of the hierarchical scheme converges to the solution of the microscopic time stepper when the macroscopic simulation interval of size $\Delta t$ decreases to zero. Key is that the simulation horizon $T$ is fixed for good comparison between different time step sizes $\delta t$ and $\Delta t$.

The proof of the consistency result (Theorem~\ref{thm:consistency}) holds under two general assumptions. These are assumptions on the solutions of the microscopic and macroscopic models, as well as the matching step. These two assumptions hold for all methods in the sections below, although we do not prove them specifically every time.

\begin{assumption} \label{ass:1}
    The solutions $w(t, x)$ and $W(t, x)$ of the microscopic~\eqref{e:moment_system} and macroscopic models, respectively, are continuous on $[0, T] \times \mathbb{R}^d$.
\end{assumption}

\begin{assumption} \label{ass:2}
    The matching operator~\eqref{eq:general_matching} is continuous in $\Delta t$, and specifically
    \begin{equation*}
    \lim_{\Delta t \to 0} f^{n+1} = f^{n, *}.
    \end{equation*}
    for any $n \in \mathbb{N}$ and $\delta t > 0$.
\end{assumption}

We can now prove the following consistency result.

\begin{theorem} \label{thm:consistency}
    Assume $T \in \mathbb{R}^{+}$ is fixed. Under Assumptions~\ref{ass:1} and~\ref{ass:2}, the solution of the hierarchical micro-macro acceleration $f_T^{\delta t, \Delta t}$ is continuous in $\Delta t$ and converges to the solution of the microscopic model $f_T^{\delta t}$ when $\Delta t$ decreases to zero.
\end{theorem}
\begin{proof}

For a fixed end time $T$, there is a maximum of $n_{max} = \left \lfloor \frac{T}{\delta t}\right \rfloor$ microscopic intervals of size $\delta t$ that fit inside $[0,T]$. The remainder $T - n_{max} \delta t < \delta t$ can be filled up with a final microscopic simulation step. 

As $\Delta t$ goes to zero, there are a finite number of values for $\Delta t$ that may cause discontinuities. The solution of the hierarchical micro-macro acceleration method is continuous in $\Delta t$ due to Assumption~\ref{ass:1}. These values are $\delta t + \Delta t = \frac{T}{n}$, for $n = 1, \dots, n_{max}$. If $\Delta t$ is larger than this value, there are $n-1$ full hierarchical micro-macro acceleration time steps and a smaller additional micro-macro step. For $\Delta t$ smaller than this value, one needs $n$ full micro-macro steps, with some additional microscopic simulation to make end time $T$.

To prove continuity of the full hierarchical micro-macro acceleration solution in terms of $\Delta t$, we need to especially check continuity at these key values. We verify continuity by investigating the left and right limit.

The left limit is reached if $\Delta t < \frac{T}{n} - \delta t$. There is a small microscopic simulation time of $\nu = T - n \left(\delta t + \Delta t \right)$ required to reach full simulation time $T$. This extra simulation is just the microscopic moment model, which is continuous by~\ref{ass:1}. This establishes left continuity.

The right limit is reached if $\Delta t > \frac{T}{n} - \delta t$, the final macroscopic simulation time interval has size $\Delta t - \nu$.
From both Assumptions~\ref{ass:1} and~\ref{ass:2}, the macroscopic simulation and matching are continuous in $\Delta t$ (and thus $\nu$). This establishes the right limit.

Put together, we have verified that the solution of the hierarchical micro-macro acceleration is continuous in $\Delta t$, and also converges to the microscopic solution as $\Delta t$ decreases to zero.
\end{proof}
\section{New methods derived from the hierarchical micro-macro acceleration}
\label{sec:3}
Using the hierarchical micro-macro acceleration presented in the previous section, we derive three new methods based on different microscopic models, that are often used for simulation of rarefied gases. The first example in Section \ref{sec:3_mMHME} considers a non-linear moment model, the second example in Section \ref{sec:4_mMHSM} uses a linear Hermite expansion and the third example in Section \ref{sec:7_mM2DHME} details the extension to multiple dimensions in physical space and velocity space.

\subsection{Example 1: micro-Macro Hyperbolic Moment Equations (mMHME)}
\label{sec:3_mMHME}
A natural choice for a hierarchical micro-macro method is to choose a refined (accurate but potentially stiff) moment model together with a coarse (inaccurate but not stiff) macroscopic model such as the Euler equations.

\subsubsection{mMHME: Micro model}
\label{sec:3_mMHME_micro}
The micro model uses the hyperbolic moment equations \cite{Cai2013b,Koellermeier2014} based on an expansion of the distribution function $f(t,x,c)$ in a Hermite series around local equilibrium \cite{Fan2016} with basis functions $\phi^{[\vel(t,x),\theta(t,x)]}_{\alpha}$
\begin{equation} \label{e:expansion}
    f(t,x,c)=\sum_{\alpha=0}^{M-1} f_{\alpha}(t,x)
    \phi^{[\vel(t,x),\theta(t,x)]}_{\alpha}\left(\frac{c-\vel(t,x)}{\sqrt{\theta(t,x)}}\right),
\end{equation}
where $f_{\alpha}(t,x)$ are the Hermite expansion coefficients and $M \in \mathbb{N}$ denotes the order of the expansion. A large $M$ results in an accurate micro model. The basis functions $\phi^{[\vel,\theta]}_{\alpha}$ are the weighted Hermite polynomial functions
\begin{equation}\label{e:grad-basisfunction}
    \phi^{[\vel,\theta]}_{\alpha}(c) = (-1)^{\alpha} \frac{\mathrm{d}^\alpha}{\mathrm{d}
        c^{\alpha}}\omega^{[\vel,\theta]}(c),
        \quad
        \omega^{[\vel,\theta]}(c)=\frac{1}{\sqrt{2\pi\theta}}\exp\left(
        -\frac{c^2}{2}\right).
\end{equation}
The ansatz is highly non-linear as the basis functions depend on the macroscopic variables $\vel(t,x)$ and $\theta(t,x)$, which are in turn moments of the distribution function \eqref{e:macroscopic_quantities}-\eqref{e:macroscopic_quantities_end}. Note that the ansatz \eqref{e:expansion} includes the macroscopic variables $\vel, \theta$. To ensure that the ansatz \eqref{e:expansion} has the correct moments corresponding to the macroscopic variables, we need to set ${f_0=\rho, f_1=0, f_2=0}$. These three variables can thus be directly inserted in the ansatz and the remaining set of micro variables reads 
\begin{equation}\label{e:vars_HME}
\var = \left(\rho, \vel, \theta, f_3, \ldots, f_{M-1}\right) \in \mathbb{R}^{M}.
\end{equation}


The evolution equations for the micro variables are of the form \eqref{e:moment_system} with system matrix $\Vect{A}_{\textrm{HME}} \in \mathbb{R}^{M\times M}$ defined by \cite{Koellermeier2017a}
\begin{equation}
\label{e:HME_A}
    \Vect{A}_{\textrm{HME}} = \setlength{\arraycolsep}{1pt} \left(
    \begin{array}{cccccccc}
        \vel & \rho &  &  &   &   &   &  \\
        \frac{\theta}{\rho} & \vel & 1 &  &  &   &   &   \\
        & 2 \theta & \vel & \frac{6}{\rho} &  &  &   &   \\
        & 4 f_3 & \frac{\rho \theta}{2} & \vel & 4 &  &  &  \\
        -\frac{\theta f_3}{\rho} & 5f_4 & \frac{3f_3}{2} & \theta & \vel & 5 &  &    \\
        \vdots & \vdots & \vdots & \vdots &  \ddots & \ddots & \ddots &  \\
        -\frac{\theta f_{M-3}}{\rho} & (M-1) f_{M-1} & \frac{\left(M-3\right)f_{M-3}+ \theta f_{M-5}}{2} & -\frac{3f_{M-4}}{\rho}  &  & \theta & \vel & M-1 \\
        -\frac{\theta f_{M-2}}{\rho} & \textcolor{myred}{0} & \textcolor{myred}{-f_{M-2}} + \frac{\theta f_{M-4}}{2} & -\frac{3f_{M-3}}{\rho} &  &  & \theta & \vel \\
    \end{array}
    \right) \setlength{\arraycolsep}{6pt},
\end{equation}
and the  right-hand side source term $\vect{S}(\var) \in \mathbb{R}^{M}$ modeling collisions using the simple BGK model \cite{Bhatnagar1954} as
\begin{equation}
\label{e:BGK_term}
    - \frac{1}{\epsilon} \vect{S}(\Var)  = - \frac{1}{\epsilon} diag \left( 0,0,0,1, \ldots,1 \right) \, \var,
\end{equation}
for relaxation time $\epsilon \in \mathbb{R}_+$. Note how the source term leads to a relaxation of the coefficients $f_i, i \geq 3,$ to zero, which is the relaxation to equilibrium, in which the distribution function $f(t,x,c)$ is in the form of a Maxwellian \eqref{e:Maxwellian} and characterized solely by the first three moments $\rho, \vel, \theta$.
We note that other hyperbolic moment models of the form \eqref{e:moment_system} can be readily used \cite{Koellermeier2015,Koellermeier2014}.

The spectral analysis in \cite{Koellermeier2021} has revealed a spectral gap for small relaxation times due to the different times scales of the transport terms and the collision terms. The fast relaxation of the higher moments requires small time steps for the stability of standard time stepping schemes while the transport terms only require a standard CFL-type time step size. 
The micro step thus uses a small time step size $\delta t$ according to the spectral properties of the micro model proved in \cite{Koellermeier2021}.

\subsubsection{mMHME: Macro model}
\label{sec:3_mMHME_macro}
The macro model uses the Euler equations of fluid dynamics, given by \eqref{decoupled-scheme}. The macro variables are then given by the equilibrium values density, velocity, and temperature, i.e., $\Var = \left(\rho, \vel, \theta\right) \in \mathbb{R}^{3}$.

The macro model uses a large CFL-type time step size $\Delta t$ according to the respective properties of the macro model, which is free from stiff, relaxing moments, compare \ref{e:BGK_term} and \cite{Koellermeier2021}.

\subsubsection{mMHME: Restriction}
\label{sec:3_mMHME_restriction}
The restriction of the intermediate micro solution to the macro solution is a simple extraction of the density, velocity, and temperature values:
\begin{equation*}
    \Var^{n,*} = (\rho^{n,*}, \vel^{n,*}, \theta^{n,*}) \in \mathbb{R}^3,
\end{equation*}
where the single entries are simply carried over from the intermediate micro solution $\var^{n,*}$.

\subsubsection{mMHME: Matching}
\label{sec:3_mMHME_matching}
To obtain a new micro solution, matching uses the consistency of moments with the macro solution and the minimum distance to the intermediate micro solution. Two choices are needed: the number of moments $L$ to be matched and the measure for the distance $d$.

The simplest choice for the number of moments $M-L$ to be matched is to match only the first three moments of the macro model, corresponding to density $\rho$, velocity $\vel$, and temperature $\theta$ from the macro Euler model. This already fixes the first three entries of the new micro solution. The others depend on the choice of the distance function $d$.

The simplest choice for the distance function between the new micro solution and the intermediate micro solution is to use a weighted $L^2$ norm, i.e.,
\begin{equation}\label{eq:L^2_norm}
    d(f, f^{n,*}) = || f^{n,*} - f ||_{2,\omega} = \int_\mathbb{R} \left(f^{n,*}(c) - f(c)\right)^2 \omega^{[\vel^{n+1},\theta^{n+1}]}(c)^{-1} \, dc.
\end{equation}

Due to the expansion in basis functions and coefficients \eqref{e:expansion}, the weighted $L^2$ distance can be explicitly computed. The intermediate micro solution $f^{n,*}$ uses basis functions $\phi^{n,*}_{i}$ with intermediate macro values $\rho^{n,*}, \vel^{n,*}, \theta^{n,*}$ and the new micro solution uses basis functions $\phi_{i}$ with the new macro values already extracted from the macro solution $\rho^{n+1}, \vel^{n+1}, \theta^{n+1}$. The distance can then be expressed as
\begin{eqnarray*}
     && \int_\mathbb{R} \left(f^{n,*}(c) - f(c)\right)^2 \omega^{[\vel^{n+1},\theta^{n+1}]}(c)^{-1} \, dc \\
    &=& \int_\mathbb{R} \left(f^{n,*}(c)^2 - 2 f^{n,*}(c) f(c) + f(c)^2\right) \omega^{[\vel^{n+1},\theta^{n+1}]}(c)^{-1} \, dc \\
    &=& \int_\mathbb{R} \left(\sum_{i,j=0}^{M-1} f^{n,*}_i f^{n,*}_j \phi^{n,*}_{i} \phi^{n,*}_{j} - 2 \sum_{i,j=0}^{M-1} f^{n,*}_i f_j \phi^{n,*}_{i} \phi_{j} + \sum_{i,j=0}^{M-1} f_i f_j \phi_{i} \phi_{j} \right) \omega^{[\vel^{n+1},\theta^{n+1}]}(c)^{-1} \, dc \\
    &=& C  - 2 \sum_{i,j=0}^{M-1} f^{n,*}_i f_j B_{i,j} + \sum_{i,j=0}^{M-1} f_i f_j A_{i,j},
\end{eqnarray*}
for constant $C:= \int_\mathbb{R} \left(\sum_{i,j=0}^{M-1} f^{n,*}_i f^{n,*}_j \phi^{n,*}_{i} \phi^{n,*}_{j} \right) \omega^{[\vel^{n+1},\theta^{n+1}]}(c)^{-1} \, dc$, and using the definitions
\begin{equation}\label{eqn:L^2_matrices}
   B_{i,j} := \int_\mathbb{R} \phi^{n,*}_{i} \phi_{j} \omega^{[\vel^{n+1},\theta^{n+1}]}(c)^{-1} \, dc, \quad A_{i,j} := \int_\mathbb{R} \phi_{i} \phi_{j} \omega^{[\vel^{n+1},\theta^{n+1}]}(c)^{-1} \, dc.
\end{equation}
The search for the minimum distance then requires the derivatives with respect to the remaining $f_i, i = 3, \ldots, M-1$, which is given by
\begin{equation} \label{eqn:L^2_derivative}
   \frac{\partial || f^{n,*} - f ||_{2,\omega}}{\partial f_i} = - 2 \sum_{j=0}^{M-1} f^{n,*}_j B_{i,j} + 2 \sum_{j=0}^{M-1} f_j A_{i,j}.
\end{equation}
Setting the derivatives \eqref{eqn:L^2_derivative} to zero, the solution coefficients $f^{n+1}$ are given by 
\begin{equation}\label{eqn:L^2_sol}
   f^{n+1} = A^{-1} B f^{n,*},
\end{equation}
where $A, B \in \mathbb{R}^{M \times M}$ are the matrices defined in \eqref{eqn:L^2_matrices}. Note that $A=A\left(u^{n+1}, \theta^{n+1}\right)$ only depends on the new macro values $\vel^{n+1}$, $\theta^{n+1}$, while $B=B\left(u^{n+1}, \theta^{n+1}, u^{n,*}, \theta^{n,*}\right)$ additionally depends on the intermediate micro values $\vel^{n,*}$, $\theta^{n,*}$. All expressions can be analytically precomputed and evaluated at runtime. In case of an orthonormal basis $\phi_{i}$, we obtain $A_{i,j} = \delta_{i,j}$, such that the matrix $A$ is the unit matrix and no inversion is necessary. In case of an orthogonal matrix, only row-wise scaling needs to be performed to invert $A$. Similarly, $B$ is typically an upper triangular matrix. The computation of the solution $f$ from \eqref{eqn:L^2_sol} is thus no more than a simple matrix vector product.

\begin{remark}
    It is possible to use a moment model as the macro model, too. For example a moment model of the same type as the micro model \eqref{e:HME_A}, but with fewer moments. This leads to a more restrictive time step size for the macro model due to the right-hand side relaxation terms, but might be beneficial for the accuracy of the model. The restriction is then simply performed using more moments of the intermediate micro solution. Due to the structure of the expansion with orthogonal Hermite functions, no moments need to be explicitly computed and the macro values can directly be extracted from the intermediate micro solution. The matching step can still be performed in the same way as described, with more coefficients extracted from the macro model.
\end{remark}

\subsection{Example 2: micro-Macro Hermite Spectral Model (mMHSM)}
\label{sec:4_mMHSM}
For small velocities $\vel$ and not too large temperatures $\theta$, linearized moment models, such as described in \cite{Fan2020a,Koellermeier2021}, can be used as a micro model within the hierarchical micro-macro acceleration. These linearized models are simpler to solve and their linear structure will also make the matching step with the macro model even simpler compared to section \ref{sec:3_mMHME}.

\subsubsection{mMHSM: Micro model}
\label{sec:4_mMHSM_micro}
Using basis functions that do not depend on the local velocity $\vel$ and temperature $\theta$, a linearized model, called the Hermite Spectral Model (HSM), is derived \cite{Koellermeier2021}. The model still uses the same basis as the non-linear HME model from section \ref{sec:3}, but can be seen as a linearization around $\vel = 0$, $\theta = 1$, for which the equilibrium Maxwellian \eqref{e:Maxwellian} is simply a Gaussian centered at $\vel = 0$ with variance $\theta = 1$.

The expansion can be written as
\begin{equation} \label{e:expansionHSM}
    f(t,x,c)=\sum_{\alpha = 0}^{M-1} f_{\alpha}(t,x) \mathcal{H}_{\alpha}(c),
\end{equation}
with weighted Hermite basis functions $\mathcal{H}_{\alpha}$ that do not depend on $\vel, \theta$ as
\begin{equation}\label{e:grad-basisfunctionHSM}
    \mathcal{H}_{\alpha}(c) = \frac{1}{\sqrt{2 \pi}} \exp \left( -\frac{c^2}{2} \right) He_{\alpha}(c) \cdot \frac{1}{\sqrt{2^{\alpha} {\alpha}!}},
\end{equation}
where $He_{\alpha}$ is the normalized Hermite polynomial of degree $\alpha$.

The constraints for the reproduction of the first three macroscopic moments are
\begin{equation} \label{e:constraints_Grad_HSM}
    f_{0} = \rho, ~~ f_{1} = \rho \vel, ~~ f_{2} = \frac{1}{\sqrt{2}} \left( \rho \theta + \rho \vel^2 - \rho \right), 
\end{equation}
so that the vector of variables can be written, e.g., as $\vect{f}= \left(f_0,\ldots,f_{M-1}\right) \mathbb{R}^M$.

The evolution equations for the HSM model can be written in the form of \eqref{e:moment_system} with system matrix $\Vect{A}_{\textrm{HSM}} \in \mathbb{R}^{M\times M}$ defined by
\begin{equation}
\label{e:grad_A_HSM}
\Vect{A}_{\textrm{HSM}} = \left(
  \begin{array}{ccccc}
     & 1 &   &   &   \\
     1 &  & \sqrt{2}  &   &   \\
     & \sqrt{2} &   &  \ddots &   \\
     &  & \ddots  &   & \sqrt{M-1}  \\
     &  &   & \sqrt{M-1}  &
  \end{array}
\right).
\end{equation}

The right-hand source term $\vect{S}\left( \Var \right) \in \mathbb{R}^{M}$ is given by  
\begin{equation}
  - \frac{1}{\epsilon} \vect{S}_{\alpha} = - \frac{1}{\epsilon} \int_{\mathbb{R}} \left( f(t,x,c) - f_M(t,x,c) \right) \psi_{\alpha}(c) \,dc, \quad \textrm{ for } \psi_{\alpha}(c) = He_{\alpha}(c) \cdot \frac{1}{\sqrt{2^{\alpha} {\alpha}!}}.
\end{equation}
Using \eqref{e:expansionHSM} the source term can be computed fully analytically. We omit the details of the derivation here for conciseness. Note that the right-hand side still exhibits the same relaxation behavior as the non-linear model from the previous section.

The micro step therefore again uses a small time step size $\delta t$ according to the spectral properties of the HSM model proved in \cite{Koellermeier2021}.

\subsubsection{mMHSM: Macro model}
\label{sec:4_mMHSM_macro}
The macro model again uses the Euler equations of fluid dynamics, given by \eqref{decoupled-scheme}. The macro variables are the equilibrium values density, velocity and, temperature, such that $\Var = \left(\rho, \vel, \theta\right) \in \mathbb{R}^{3}$.

\subsubsection{mMHSM: Restriction}
\label{sec:4_mMHSM_restriction}
The restriction of the intermediate micro solution to the macro solution is again a simple extraction of the density, velocity and temperature values:
\begin{equation*}
    \Var^{n,*} = (\rho^{n,*}, \vel^{n,*}, \theta^{n,*}) \in \mathbb{R}^3,
\end{equation*}
where the single entries are computed according to the inverted consistency constraints \eqref{e:constraints_Grad_HSM} from $\var^{n,*}$ as
\begin{equation*}
    \rho^{n,*} = f_0^{n,*}, \quad \vel^{n,*} = \frac{f_1^{n,*}}{f_0^{n,*}}, \quad \theta^{n,*} = 1 + \sqrt{2} \frac{f_2^{n,*}}{f_0^{n,*}}  - \left(\frac{f_1^{n,*}}{f_0^{n,*}}\right)^2.
\end{equation*}

\subsubsection{mMHSM: Matching}
\label{sec:4_mMHSM_matching}
The matching step can be written in the same form as for the HME model in section \ref{sec:3_mMHME_matching}, again matching the first three moments and minimizing the weighted $L^2$ distance. 

The first three coefficients $f_0^{n+1}, f_1^{n+1}, f_2^{n+1}$ of the new micro solution are computed by matching with the exact moments of the macro model using the consistency constraints \eqref{e:constraints_Grad_HSM}, i.e.,
\begin{equation} \label{e:constraints_Grad_HSM_matching}
    f_{0}^{n+1} = \rho^{n+1}, ~~ f_{1}^{n+1} = \rho^{n+1} \vel^{n+1}, ~~ f_{2}^{n+1} = \frac{1}{\sqrt{2}} \left( \rho^{n+1} \theta^{n+1} + \rho^{n+1} (\vel^{n+1})^2 - \rho^{n+1} \right). 
\end{equation}

The rest of the coefficients is computed via minimizing the distance to the intermediate micro solution. However, both the basis functions as well as the weight function do not depend on $\vel,\theta$ in this linearized setup, which simplifies the matching. The linear system for the new micro coefficients $f^{n+1}$ reads
\begin{equation*} 
    \sum_{j=0}^{M-1} f^{n,*}_j B_{i,j} = \sum_{j=0}^{M-1} f_j^{n+1} A_{i,j}, i = 3, \ldots, M-1,
\end{equation*}
with
\begin{equation}\label{eqn:L^2_matrices_HSM}
   B_{i,j} := \int_\mathbb{R} \phi^{n,*}_{i} \phi_{j} \omega(c)^{-1} \, dc = \delta_{i,j}, \quad A_{i,j} := \int_\mathbb{R} \phi_{i} \phi_{j} \omega(c)^{-1} \, dc = \delta_{i,j}.
\end{equation}
The solution is simply
\begin{equation}\label{eqn:L^2_sol_HSM}
   f_i^{n+1} = f^{n,*}_i, \quad i = 3, \ldots, M-1,
\end{equation}
which means that the remaining coefficients of the intermediate micro solution are simply carried over to the next micro step.

\subsection{Example 3: micro-Macro 2D Hyperbolic Moment Equations (mM2DHME)}
\label{sec:7_mM2DHME}
In this section, we extend the 1D mMHME method derived in \ref{sec:3_mMHME} to the multi-dimensional case. We exemplarily consider the 2D case as the extension to higher dimensions is straightforward. While the hierarchical micro-macro acceleration stays the same, the method slightly differs due to the multi-dimensional nature.

\subsubsection{mM2DHME: Micro model}
The 2D micro model uses 2D hyperbolic moment equations, see e.g. \cite{Koellermeier2017d}, based on an expansion of the distribution function $f(t,\vect{x},\vect{c})$ in a Hermite series around local equilibrium with basis functions $\phi^{[\vect{\vel}(t,\vect{x}),\theta(t,\vect{x})]}_{\vect{\alpha}}$
\begin{equation} \label{e:2Dexpansion}
    f(t,\vect{x},\vect{c})=\sum_{\vect{\alpha}\leq M}^{} f_{\vect{\alpha}}(t,\vect{x})
    \phi^{[\vect{\vel}(t,\vect{x}),\theta(t,\vect{x})]}_{\vect{\alpha}}\left(\vect{c}\right),
\end{equation}
where $\vect{\alpha}\in \mathbb{N}^2$ is a multi-index, $f_{\vect{\alpha}}(t,\vect{x})$ are the Hermite expansion coefficients, and $M \in \mathbb{N}$ denotes the order of the expansion, where full moment theories are used. The basis functions $\phi^{[\vect{\vel},\theta]}_{\alpha}$ are the tensor products of weighted Hermite polynomial functions
\begin{equation}\label{e:2Dgrad-basisfunction}
    \phi^{[\vect{\vel},\theta]}_{\vect{\alpha}}(c) =  \prod_{i=1}^2 \frac{1}{\sqrt{2\pi\theta}}\exp\left(
        -\frac{(\vect{c}-\vect{\vel})^2}{2\theta}\right)\cdot \left(2\theta\right)^{-\frac{\alpha_i}{2}} He_{\alpha_i}\left(\frac{c_i-\vel_i}{\sqrt{\theta}}\right),
\end{equation}
where $He_{\alpha_i}(c)$ is the standard Hermite polynomial of degree $i$.

The ansatz is highly nonlinear due to the dependence on $\vect{\vel}, \theta$. In 2D, the compatibility conditions, which ensure that the correct mass, momentum and energy are recoverd, read: 
$$f_{0,0} = \rho, \quad f_{1,0} = f_{0,1} = 0, \quad f_{2,0} + f_{0,2} = 0.$$
We focus on the case $M=3$, where the full moment ansatz results in the following ten variable vector, which was used, e.g., in the simulations in \cite{Koellermeier2021,Koellermeier2018b}. 
\begin{equation}
\label{e:2Dvars_full}
    \var = \left( \rho, \vel_x, \vel_y, \frac{p_1}{2}, f_{1,1}, \frac{p_2}{2}, f_{3,0}, f_{2,1}, f_{1,2}, f_{0,3}\right)^T.
\end{equation}
for $\frac{p_1}{2} = \frac{\rho \theta}{2}+f_{2,0}$, $\frac{p_2}{2} = \frac{\rho \theta}{2}+f_{0,2}$, and $f_{i,j} = f_{i\vect{e}_1+j\vect{e}_2}$.

Inserting the ansatz \eqref{e:2Dexpansion} into the 2D version of the kinetic equation \eqref{e:BTE} and subsequently projecting onto orthogonal basis functions yields the closed moment equations. The hyperbolic version of the moment equations called Hyperbolic Moment Equations (HME), can be written as  
\begin{equation}
\label{e:2Dmoment_system}
    \frac{\partial \var}{\partial t}  + \Vect{A}_{x} \frac{\partial \var}{\partial {x}} + \Vect{A}_{y} \frac{\partial \var}{\partial {y}} = \Vect{S}(\var),
\end{equation}
where $\Vect{S}(\var)$ is the right-hand side source term. For conciseness, the matrices $\Vect{A}_{x}$, $\Vect{A}_{x}$ and the source term $\Vect{S}(\var)$ are omitted here and can be found in \cite{Koellermeier2017d,Koellermeier2018b,Koellermeier2021} with details and alternative 2D moment models.

Similarly to the 1D model, the 2D model is characterized by a spectral gap for small relaxation times, due to the different scales of the transport and collision terms. The micro model therefore requires a small step size $\delta t$ for stability as proved in \cite{Koellermeier2021}.

\subsubsection{mM2DHME: Macro model}
As macro model, we use the 2D Euler equations, which include the following variables
\begin{equation}
\label{e:Eulervars_full}
    \Var = \left( \rho, \vel_x, \vel_y, p\right)^T,
\end{equation}
and the model equations are given by 
\begin{equation}\label{e:2DEuler}
        \partial_t \Var +
      \left( \begin{array}{cccc}
        \vel_x & \rho & 0 & 0\\
        0 & \vel_x & 0 & \frac{1}{\rho}\\
        0 & 0 & \vel_x & 0\\
        0 & \gamma p & 0 & \vel_x \\
      \end{array} \right)
       \partial_x \Var +
      \left( \begin{array}{cccc}
        \vel_y & 0 & \rho & 0\\
        0 & \vel_y & 0 & 0\\
        0 & 0 & \vel_y & \frac{1}{\rho}\\
        0 & 0 & \gamma p & \vel_y \\
      \end{array} \right)
       \partial_y \Var = \left( \begin{array}{c}
        0 \\
        0 \\
        0 \\
        0 \\
      \end{array} \right),
\end{equation}
where $\gamma = 2$ for 2D ideal gases. Note that the collision term vanishes for the Euler equation, which means that the fast scales are eliminated from the model. Therefore, the Euler equations can use a macroscopic time step size $\Delta t$ that is only limited by a standard CFL condition or accuracy considerations.

\subsubsection{mM2DHME: Restriction}
\label{sec:7_mM2DHSM_restriction}
The restriction of the intermediate micro solution to the macro solution only requires setting the macro solution
\begin{equation*}
    \Var^{n,*} = (\rho^{n,*}, \vel_x^{n,*}, \vel_y^{n,*}, p^{n,*}) \in \mathbb{R}^4,
\end{equation*}
where the macro pressure is computed from the components of the micro solution as
\begin{equation*}
    p^{n,*} = \frac{p_1}{2} + \frac{p_2}{2}.
\end{equation*}

\subsubsection{mM2DHME: Matching}
During matching, the coefficients of the expansion \eqref{e:2Dexpansion} are computed based on the macroscopic moments obtained from the Euler solution $\Var^{n}$ and the intermediate (prior) micro solution $\var^{n,*}$. We first compute the temperature $\theta^n$ used in the basis functions of the expansion \eqref{e:2Dexpansion} as $\theta^n = p^n / \rho^n$ and the consistent intermediate micro temperature is $\theta^{n,*} = \left(\frac{p_1}{2}^{n,*} + \frac{p_2}{2}^{n,*}\right) / \rho^{n,*}$.

Similar to the 1D case, We can then solve the matching problem
\begin{equation}\label{eqn:2DL^2_sol}
   f^{n+1} = A^{-1} B f^{n,*}.
\end{equation}
In \eqref{eqn:2DL^2_sol}, the known intermediate (prior) coefficients $f^{n,*}$ are given by   
\begin{eqnarray*}
     && \left( f_{0,0}, f_{1,0}, f_{0,1}, f_{2,0}, f_{1,1}, f_{0,2}, f_{3,0}, f_{2,1}, f_{1,2}, f_{0,3}\right)^{n,*} \\
    &=& \left( \rho, 0, 0, \frac{1}{2}\left(\frac{p_1}{2}-\frac{p_2}{2}\right), f_{1,1}, \frac{1}{2}\left(\frac{p_2}{2}-\frac{p_1}{2}\right), f_{3,0}, f_{2,1}, f_{1,2}, f_{0,3}\right)^{n,*}.
\end{eqnarray*}
The matrices $A$ and $B$ in \eqref{eqn:2DL^2_sol} are precomputed as
\begin{equation}\label{eqn:2DL^2_matrices}
   B_{i,j} := \iint_{\mathbb{R}^2} \phi^{n,*}_{i} \phi_{j} \omega^{n+1}(\vect{c})^{-1} \, \vect{dc}, \quad A_{i,j} := \iint_{\mathbb{R}^2} \phi_{i} \phi_{j} \omega^{n+1}(\vect{c})^{-1} \, \vect{dc}.
\end{equation}

The vector of micro variables $\var^{n+1}$ can be obtained thereafter using    
\begin{eqnarray*}
     && \left( \rho, \vel_x, \vel_y, \frac{p_1}{2}, f_{1,1}, \frac{p_2}{2}, f_{3,0}, f_{2,1}, f_{1,2}, f_{0,3}\right)^{n+1} \\
    &=& \left( \rho, \vel_x, \vel_y, \frac{\rho \theta}{2} + f_{2,0}, f_{1,1}, \frac{\rho \theta}{2} + f_{0,2}, f_{3,0}, f_{2,1}, f_{1,2}, f_{0,3}\right)^{n+1}.
\end{eqnarray*}

The solution of the matching problem can be implemented analytically such that a small matrix vector multiplication is sufficient to compute the new coefficients. An extension for higher dimensions or larger moment models $M>3$ is straightforward.

\section{Existing methods compared to hierarchical micro-macro acceleration}
\label{sec:3b}
In this section, we compare two existing methods to the hierarchical micro-macro acceleration. For conciseness, we consider the 1D case again.

\subsection{Example 4: Projective Integration (PI)}
\label{sec:5_PI}
Instead of chosing a physical macro model, it is also possible to only extrapolate the micro solution over a large time step. In the literature, this is known as Projective Integration (PI)~\cite{Gear2003, Gear2003e}. In this section, we describe PI in the setting of the hierarchical micro-macro acceleration for comparison. Note that many adaptations to the PI exist and we cover the standard version, see e.g. \cite{Koellermeier2021,Lafitte2017,Melis2019,Melis2016}.

\subsubsection{PI: Micro model}
\label{sec:5_PI_micro}
Both the HME or the HSM model can serve as micro model for the PI method written as a hierarchical micro-macro acceleration. See \cite{Koellermeier2021} for a comparison of both moment models when using PI. In the following, we assume that the HME method is used because it has the slightly simpler consistency constraints ${f_0=1, f_1=0, f_2=0}$, such that the vector of variables reads 
$$\var = \left(\rho, \vel, \theta, f_3, \ldots, f_{M-1}\right) \in \mathbb{R}^{M}.$$
Using this micro model, we compute from the old micro solution $\var^n$ the intermediate micro solution $\var^{n,*}$ by applying the micro model for a number of small time steps with time step size $\delta t$.

Note that for stability several subsequent time steps of size $\delta t$ of the micro model are necessary, so that the fast modes are sufficiently damped \cite{Koellermeier2021}. We can therefore denote $\var^{n,k}$ as the subsequent intermediate solutions for $k = 1, \ldots, K$, with $\var^{n,0}=\var^{n}$ and $\var^{n,*} = \var^{n,K}$. The choice of the number of micro steps $K$ is typically made according to the stability properties, i.e. the spectral gap of the micro model. A larger value of $K$ makes the method more stable as the fast modes of the micro model are damped further before the macro step, see \cite{Koellermeier2021}.

\subsubsection{PI: Macro model}
\label{sec:5_PI_macro}
As macro model, we do not use a physical model that evolves the non-stiff modes. Instead, we assume that the stiff modes have been sufficiently damped by the (one or several) applications of the micro model. The evolution of the macro model is then a simple extrapolation (or projection) of the micro solution over the remainder of a chosen larger time step $\Delta t$. Assuming $K$ small time steps the projection step can be written as
\begin{equation}\label{eqn:PI_proj_K}
   \Var^{n+1} = \Var^{n} + (\Delta t - K \delta t) \frac{\var^{n,K} - \var^{n,K-1}}{\delta t}.
\end{equation}

Note that the extrapolation leads to negligible overhead and is much faster than the other models using one time step of an Euler solver or more complex models.

\subsubsection{PI: Restriction}
\label{sec:5_PI_restriction}
As the micro variables are only extrapolated over the remainder of the large time step during the macro step, the restriction is the identity operator, i.e.,
\begin{equation}\label{eqn:PI_restriction}
   \Var^{n,*} = \var^{n,*}.
\end{equation}

\subsubsection{PI: Matching}
\label{sec:5_PI_matching}
Similarly to the restriction operator, matching uses the identity operator:
\begin{equation}\label{eqn:PI_matching}
   \var^{n+1} = \Var^{n+1}.
\end{equation}

Especially from the restriction and matching operators it is clear that the PI method is a very simple method that can be implemented with almost negligible overhead in comparison to a simple micro solver. As the application of the macro step does not use physical information, however, physical accuracy could be lost during this extrapolation step. Furthermore, PI requires a careful prior investigation of the micro model's stability properties \cite{Melis2019,Lafitte2017,Melis2016}. For the moment models discussed here, this has been done in \cite{Koellermeier2021} for simple right-hand side collision operators. Based on the stability analysis of the micro model, more precisely the positions of the fast and slow eigenvalues, the number of micro steps $K$ needs to be chosen to ensure stability of the combined scheme.

\subsection{Example 5: Coarse Projective Integration (CPI)}
\label{sec:6}
Due to the four step procedure of the hierarchical micro-macro acceleration described in section \ref{sec:2}, many existing methods can be compared to it and small modification of one of the four steps already result in different methods. In this section, we compare the Coarse Projective Integration method (CPI)~\cite{Debrabrant2017, Gear2003e} to the hierarchical micro-macro acceleration. The CPI method can be seen in between the hierarchical micro-macro acceleration from sections \ref{sec:3_mMHME},\ref{sec:4_mMHSM} and the PI method from section \ref{sec:5_PI}.

\subsubsection{CPI: Micro model}
\label{sec:6_CPI_micro}
Again, both the HME or the HSM model can serve as micro model for the CPI method and we assume that the HME method is used for conciseness. This means that the micro variables are given by 
$$\var = \left(\rho, \vel, \theta, f_3, \ldots, f_{M-1}\right) \in \mathbb{R}^{M} $$
and the intermediate micro solution $\var^{n,*}$ is computed from the old micro solution $\var^n$ using a number of time steps with small time step size $\delta t$. We again denote the intermediate values by $\var^{n,k}$, see section \ref{sec:5_PI_micro}.

\subsubsection{CPI: Macro model}
\label{sec:6_CPI_macro}
Similarly to the PI method from section \ref{sec:5_PI}, the macro model is evolved using extrapolation. But the CPI method uses a reduced set of variables $\Var \in \mathbb{R}^L$ with $L \leq M$. 
The new macro values are then extrapolated over the remainder of the larger time step $\Delta t$ in the same way as for the PI method, see \eqref{eqn:PI_proj_K}, where we assume $K$ micro steps have been performed,
\begin{equation}\label{eqn:CPI_proj_K}
   \Var^{n+1} = \Var^{n} + (\Delta t - K \delta t) \frac{\Var^{n,K} - \Var^{n,K-1}}{\delta t}.
\end{equation}

In order to make the dimensions compatible, note that $\Var^{n,K-1} \in \mathbb{R}^L$ are now the restricted micro solutions.

\subsubsection{CPI: Restriction}
\label{sec:6_CPI_restriction}
As only a subset of the micro variables are extrapolated in the CPI method, the restriction operator is effectively a cut-off of the micro solution
\begin{equation}\label{eqn:CPI_restriction}
    \Var^{n,*} = \var^{n,*}_{1:L},
\end{equation}
so that the first $L$ moments of the intermediate micro solution are extrapolated and the remaining $M-L$ moments are not changed during the macro step. For $L=M$ the CPI method degenerates to the PI method.

\subsubsection{CPI: Matching}
\label{sec:6_CPI_matching}
As only a subset of the micro variables are extrapolated during the macro step, the new micro solution needs to be constructed according to the same matching procedure as for the HME or HSM method, respectively. This means that the first $L$ moments are carried over from the macro solution as
\begin{equation}\label{eqn:CPI_matching}
    \var^{n+1}_i = \Var^{n+1}_i, \quad i=0, \ldots, L-1,
\end{equation}
and the remaining variables of $\var^{n+1} \in \mathbb{R}^M$ are computed by minimising the distance to the intermediate micro solution $\var^{n,*}$. For the solution of the minimisation problem, we refer to \ref{sec:3_mMHME_matching}. 

In comparison with the PI method, the CPI method uses fewer variables for the macro extrapolation step but has to perform the matching step as a result. The benefit of the CPI method is that the macro step can have a less severe stability constraint and the method therefore requires less iterations of the micro model to damp out the fast modes.     
\section{On micro-macro matching}
\label{sec:4}

\subsection{Efficient implementation of matching}
\label{sec:efficient_matching}
The exact solution of the matching problem described in section \ref{sec:3_mMHME_matching} requires only a matrix vector product with an upper triangular matrix $B \in \mathbb{R}^{M \times M}$, defined in \eqref{eqn:L^2_matrices}. However, the practical evaluation can be costly for large $M$. The matrix $B$ corrects for the change of basis functions from the intermediate micro solution $\phi^{n,*}_{i}$ to the new micro solution. The new micro solution uses the same first moments as the macro solution to form its basis functions $\phi^{n+1}_{i}$. In section \ref{sec:4_mMHSM_matching}, we have seen that the matching problem is far easier to solve if the basis functions do not depend on the moments $\rho, \vel, \theta$ and the matching step simply selects the remaining coefficients from the intermediate micro solution $f^{n,*}$.

The difficult matching problem in section \ref{sec:3_mMHME_matching} can be reformulated with the help of a basis transformation that first transforms the expansion of the intermediate micro solution $f^{n,*}$ to the basis of the new micro solution $f^{n+1}$. The subsequent matching of the higher coefficients is then trivial. We will outline this step below for the 1D case.

Assuming the distribution function of the intermediate micro solution is expanded as 
\begin{equation} \label{e:expansion_intermediate}
    f^{n,*}(t,x,c)=\sum_{\alpha=0}^{M-1} f^{n,*}_{\alpha}(t,x)
    \phi^{[\vel^{n,*},\theta^{n,*}]}_{\alpha}\left(\frac{c-\vel^{n,*}}{\theta^{n,*}}\right),
\end{equation}
and the distribution function of the new micro solution is expanded as
\begin{equation} \label{e:expansion_new}
    f^{n+1}(t,x,c)=\sum_{\alpha=0}^{M-1} f^{n+1}_{\alpha}(t,x)
    \phi^{[\vel^{n+1},\theta^{n+1}]}_{\alpha}\left(\frac{c-\vel^{n+1}}{\theta^{n+1}}\right),
\end{equation}
we observe that the two representations differ in the values of the expansion coefficients $f_{\alpha}$ as well as the basis functions $\phi^{[\vel,\theta]}_{\alpha}$, with the latter depending on the two moments $\vel$ and $\theta$. Now we perform a basis transformation of the known intermediate micro solution to the basis of the new micro solution as follows
\begin{equation} \label{e:expansion_intermediate_trafo}
    f^{n,*}(t,x,c)=\sum_{\alpha=0}^{M-1} \widetilde{f^{n,*}_{\alpha}}(t,x) \phi^{[\vel^{n+1},\theta^{n+1}]}_{\alpha}\left(\frac{c-\vel^{n+1}}{\theta^{n+1}}\right),
\end{equation}
where the transformed expansion coefficients $\widetilde{f^{n,*}_{\alpha}}$ need to be found. This basis transformation is generally costly but can be computed efficiently up to numerical precision as a simple ODE solve as noted in \cite{Cai2010}. 

The weighted $L^2$ distance $d(f^{n+1}, f^{n,*})$ between the new micro solution $f^{n+1}$ and the intermediate micro solution $f^{n,*}$ is then computed as
\begin{eqnarray*}
    && d(f^{n+1}, f^{n,*}) = || f^{n,*} - f^{n+1} ||_{2,\omega} \\
    &=& \int_\mathbb{R} \left(f^{n,*}(c) - f^{n+1}\right)^2 \omega^{[\vel^{n+1},\theta^{n+1}]}(c)^{-1} \, dc \\
    &=& \int_\mathbb{R} \left[\sum_{i=0}^{M-1} \left(\widetilde{f^{n,*}_i} - f^{n+1}_i \right) \phi^{n+1}_{i} \right]^2 \omega^{[\vel^{n+1},\theta^{n+1}]}(c)^{-1} \, dc \\
    &=& \int_\mathbb{R} \sum_{i=0}^{M-1} \left(\widetilde{f^{n,*}_i} - f^{n+1}_i \right)^2 \left(\phi^{n+1}_{i}\right)^2 \omega^{[\vel^{n+1},\theta^{n+1}]}(c)^{-1} \, dc \\
    &=& \sum_{i=0}^{M-1} \left(\widetilde{f^{n,*}_i} - f^{n+1}_i \right)^2,
\end{eqnarray*}
where we used the orthonormality of the basis $\phi^{n+1}_{i}$, i.e., $\int_\mathbb{R} \left(\phi^{n+1}_{i}\right)^2 \omega^{[\vel^{n+1},\theta^{n+1}]}(c)^{-1} \, dc = 1$ in the last step.

The minimizer is thus simply given by $f^{n+1}_i = \widetilde{f^{n,*}_i}$, for $i = 3, \ldots, M-1$, i.e., the solution of the matching problem is that the remaining new micro coefficients are the corresponding transformed intermediate micro coefficients. This means that the transformation \eqref{e:expansion_intermediate_trafo} is enough to perform the matching as the new solution $f^{n+1}_i$ simply takes the transformed values $\widetilde{f^{n,*}_i}$.

Note that the solution $f^{n+1}_i = \widetilde{f^{n,*}_i}$, for $i = 3, \ldots, M-1$, does not mean that the weighted $L^2$ distance $d(f^{n+1}, f^{n,*})$ is zero, as the first three coefficients still differ due to the incompatibility of the intermediate micro solution with the consistency conditions of the new micro solution, i.e., ${f_0^{n+1}=1, f_1^{n+1}=0, f_2^{n+1}=0}$ but ${\widetilde{f_0^{n,*}} \neq 1, \widetilde{f_1^{n,*}} \neq 0, \widetilde{f_2^{n,*}} \neq 0}$ in general.

\subsection{Alternative distance measures for matching}
So far, we have only used a weighted $L^2$ distance metric to match the new microscopic moments based on the previous microscopic step. Although the $L^2$ distance is a powerful technique for deriving the matched microscopic moments analytically, it is by no means the only (pseudo) metric that can be used.

One alternative metric that has been studied before~\cite{Debrabrant2017,Vandecasteele2020} is the Kullback-Leibler divergence. This pseudo distance measures how much information some density function $f$ contains over the prior density $f^{n,*}$. The divergence is defined as 
\begin{equation*}
    KL(f, f^{n,*}) = \int_{\mathbb{R}} f(x, c) \ln\left(\frac{f(x, c)}{f^{n,*}(x,c) } \right) dc.
\end{equation*}
Matching in Kullback-Leibler divergence can then be written, in analogy to~\eqref{eq:general_matching}, as
\begin{equation}
f^{n+1} = \underset{f  \in  V(\Var^{n+1})}{\argmin} \  KL(f, f^{n,*}).
\end{equation}
This method of matching has already been used to match a set of Monte Carlo particles in the context stochastic differential equations \cite{Vandecasteele2020}.

The main drawback of this method over $L^2$ matching is that there is no explicit formula for the matched microscopic moments $\var^{n+1}$. As a consequence, we require a numerical procedure to determine these microscopic moments, making it a much more expensive alternative~\cite{Debrabrant2017}. We refer to \cite{Vandecasteele2020} for details.

\begin{remark}
Many other pseudo metrics can be used to derive a matching operator. We mention here the general R{\'e}nyi divergence of which the Kullback-Leibler divergence is a special case~\cite{VanErven2014}. Also, the Monge-Kantorovich or Wasserstein distances from optimal transport~\cite{Villani2009} are possible alternatives for matching. The reason we have mentioned the Kullback-Leibler divergence is that it is the only alternative to $L^2$ matching that has been used before in the context of micro-macro methods \cite{Vandecasteele2020}.
\end{remark}

\subsection{Numerical investigation of matching}
The necessity of matching is a result of reducing the number of variables from the micro to the macro model. After the macro step, the unknown information of the micro variables needs to be recovered. In the hierarchical micro-macro acceleration presented in this paper, this is done via minimizing the distance to the prior intermediate micro distribution function. The reduction of variables obviously reduces the complexity of the model. We now want to investigate how much accuracy is lost by using less or more variables for the matching step itself. We therefore assume an existing intermediate micro solution (also called the prior) $w^{n,*} \in \mathbb{R}^8$ of the HME model in \eqref{e:vars_HME} based on a typical bimodal distribution function given by 
\begin{equation} \label{e:bimodal}
    \rho=1, \vel=1, \theta=1, f_3=-0.2, f_4=0.1, f_5=-0.01, f_6=0.001, f_7=-0.0005.
\end{equation}

We further assume that the exact solution is a simple scaling of the prior solution, i.e. $w^{n+1} = p \cdot w^{n,*}$, with $p=1.2$ in this test case. For the matching, we consider the $L^2$ distance function and first assume that the macro model contains $L=3$ variables, such that only the $L=3$ macro moments $\rho, \vel, \theta$ are taken directly from the exact solution and all remaining $M-L=5$ non-equilibrium variables $f_{L-1}=f_2, \ldots, f_{M-1} = f_7$ are computed using matching with the prior solution $w^{n,*}$.  

\begin{figure}[htb!]
    \centering
    \begin{subfigures}
    \subfloat[Matching all non-equilibrium variables, macro model $L=3$. \label{fig:matchingSingle}
    ]{\begin{overpic}[width=0.5\linewidth]{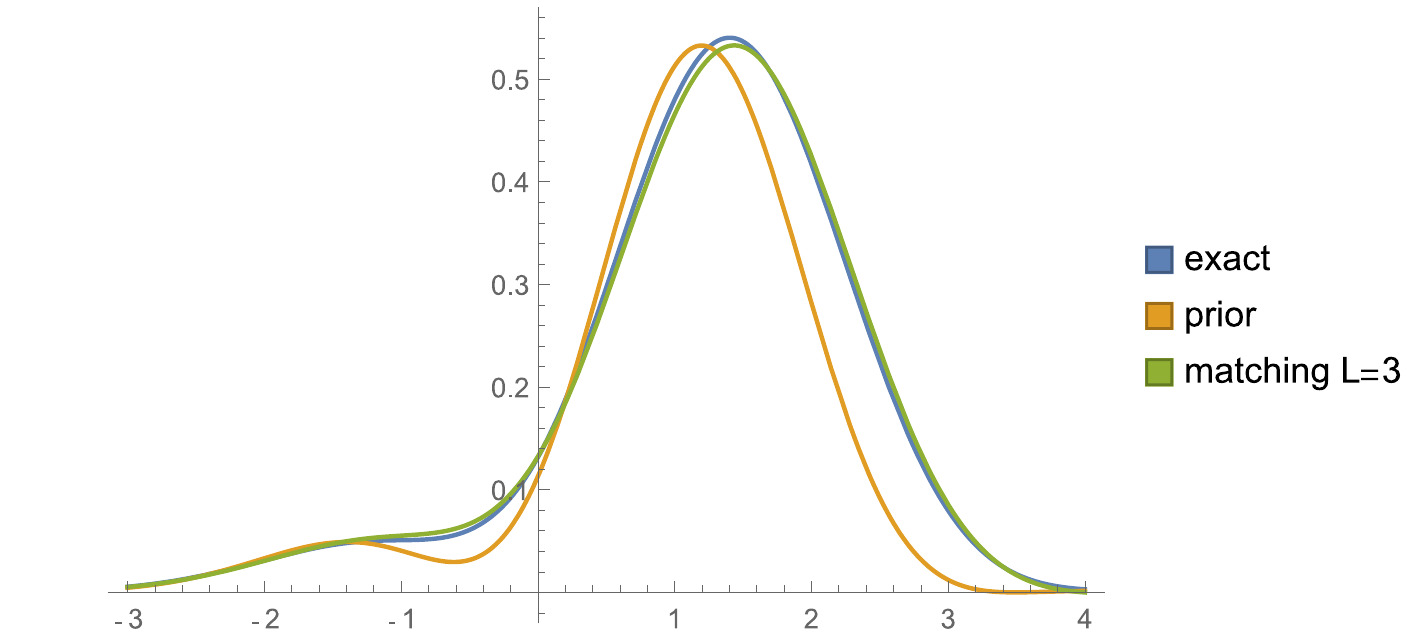}
        \put(33,43){$f$}
        \put(80,1){$c$}
    \end{overpic}}
    \subfloat[Matching different numbers of non-equilibrium variables, macro model  $L=3,\ldots,7$. \label{fig:matchingAll}
    ]{\begin{overpic}[width=0.5\linewidth]{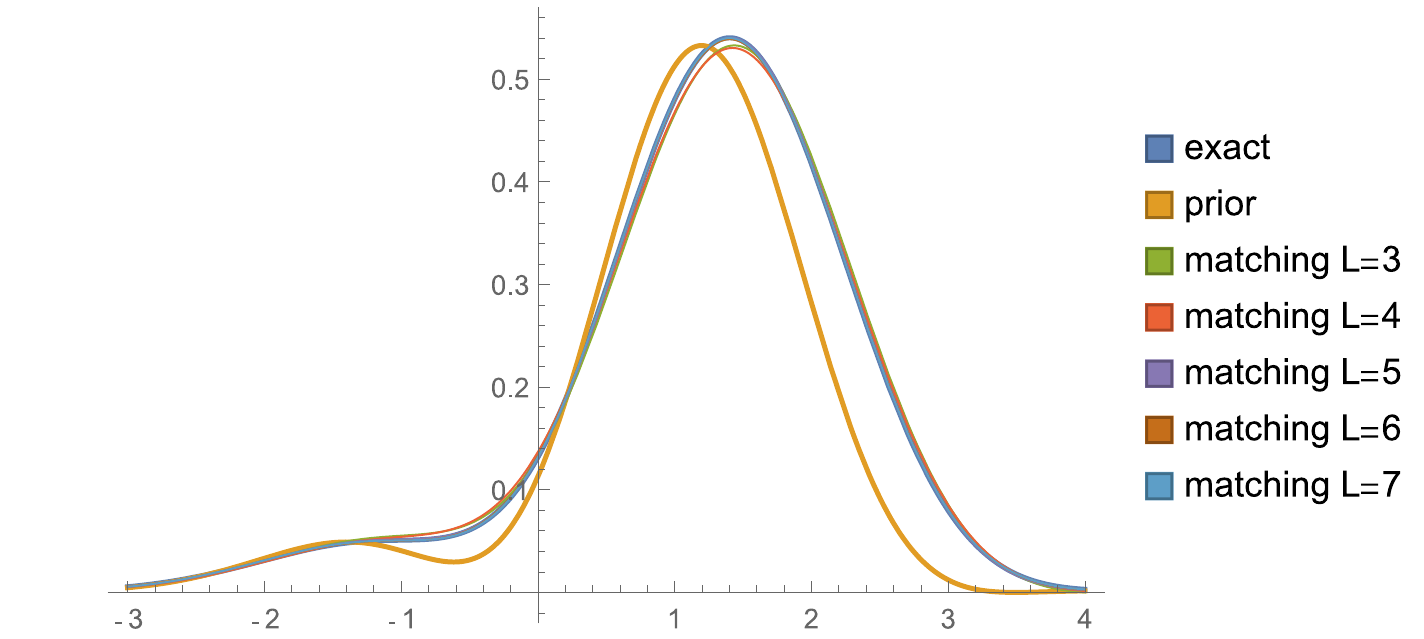}
        \put(33,43){$f$}
        \put(80,1){$c$}
    \end{overpic}}
    \end{subfigures}
    \caption{Matching all non-equilibrium variables already leads to accurate reconstruction of exact solution, while matching less variables increases accuracy at expense of more expensive macro model.}
    \label{fig:matching}
\end{figure}

\begin{figure}[htb!]
    \centering
    \begin{subfigures}
    \subfloat[Zoom in to $c\in (1,1.2)$. \label{fig:matchingZoom1}
    ]{\begin{overpic}[width=0.5\linewidth]{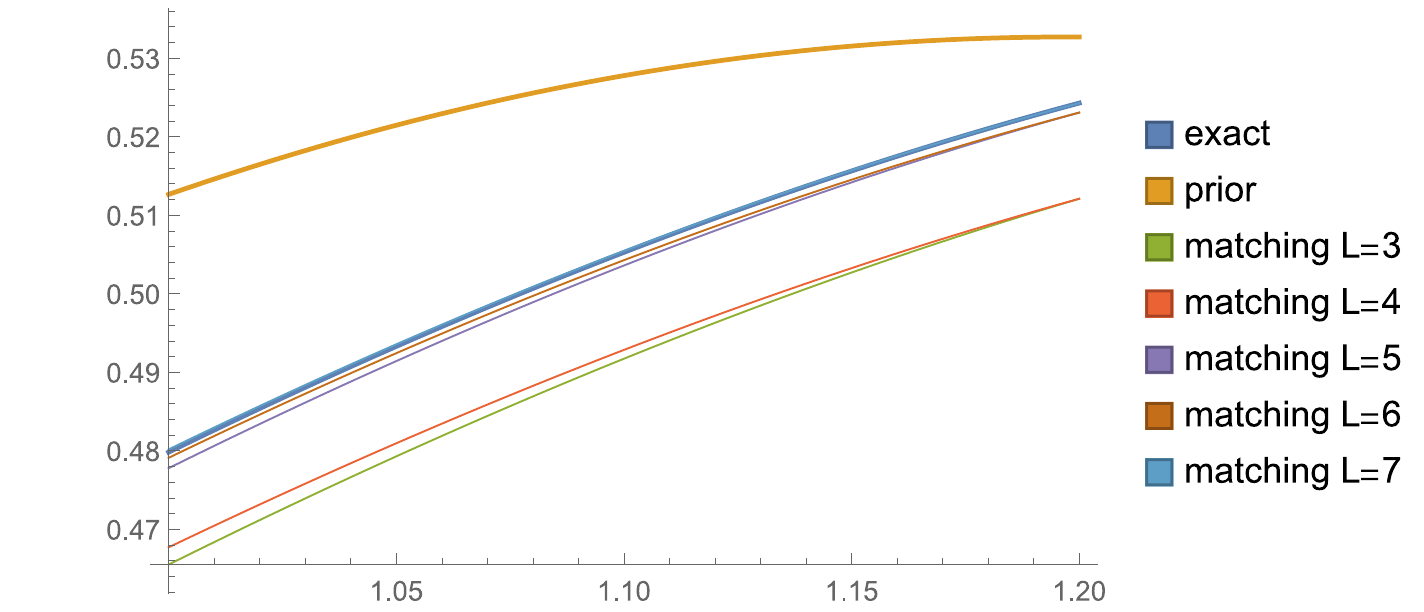}
        \put(7,43){$f$}
        \put(80,1){$c$}
    \end{overpic}}
    \subfloat[Zoom in to $c\in (0.15,0.25)$. \label{fig:matchingZoom2}
    ]{\begin{overpic}[width=0.5\linewidth]{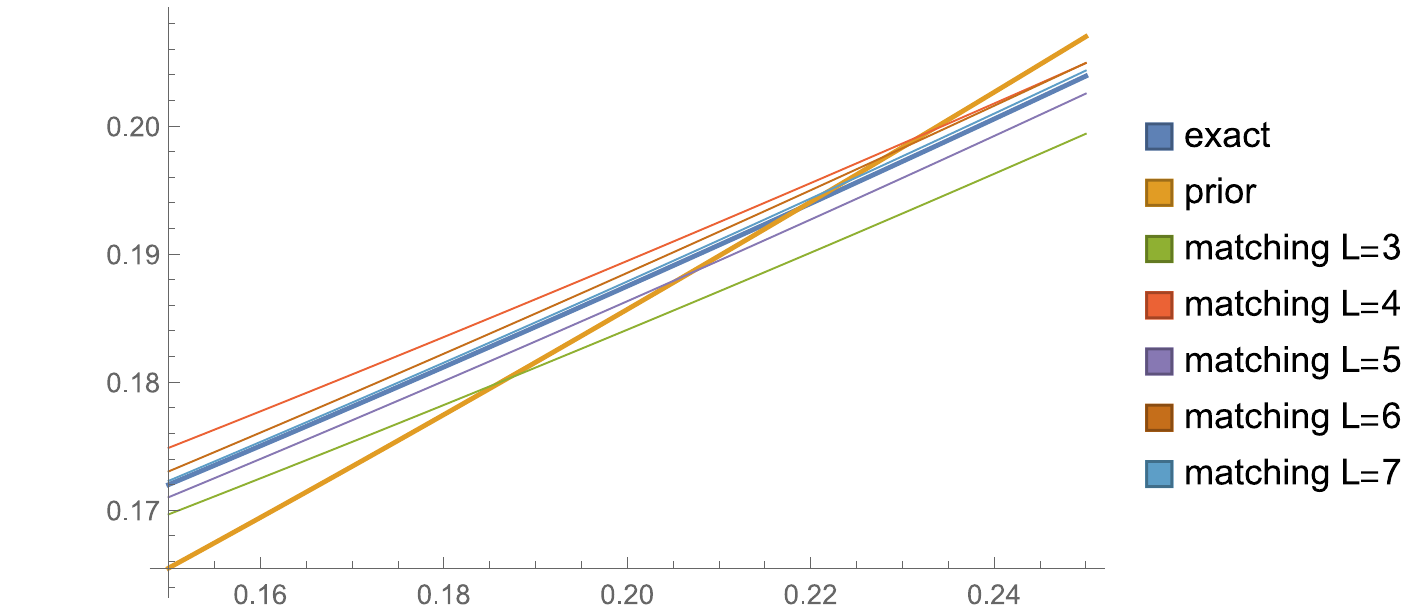}
        \put(7,43){$f$}
        \put(80,1){$c$}
    \end{overpic}}
    \end{subfigures}
    \caption{Zoom in matching different numbers of non-equilibrium variables, $L=3,\ldots,7$.}
    \label{fig:matchingZoom}
\end{figure}

Figure \ref{fig:matchingSingle} clearly shows that even though only very little information is taken from the exact $f$, the distribution function is approximated with very good quality. Matching fewer moments increases the approximation quality of the matching, as shown in figure \ref{fig:matchingAll} and the zoomed in views in figure \ref{fig:matchingZoom}. It is clear that a small number of matched variables, i.e., a large number $L$ of moments taken directly from the existing macro solution, leads to a more accurate approximation of the distribution function. However, we note that the computation of the macro solution is significantly more expensive and the reduction of moments is the reason of the hierarchical micro-macro acceleration in the first place. 

The investigation in this section shows that even when matching a large number of variables, the distribution function can still be approximated with good quality. This was also confirmed in test cases with other distribution functions than \eqref{e:bimodal}.

\subsection{Computational complexity}
\label{sec:Complexity}
The main advantage of hierarchical micro-macro acceleration is its computational speedup with respect to standard schemes close to equilibrium. In this section, we derive the typical computational complexity to allow for a numerical comparison in the next section. As all runtimes scale with the number of grid cells, we consider the complexity per grid cell.

For the reference method, we consider a micro model using $M$ moments, which corresponds to an approximate complexity of $\mathcal{O}(M^2)$ per cell and time step. For a stiff right-hand side with small relaxation time $\epsilon \ll 1$, a standard Forward Euler time stepping scheme results in a severe time step constraint $\Delta t = \mathcal{O}\left(\epsilon\right)$ leading to $\mathcal{O}\left(\frac{1}{\epsilon}\right)$ time steps per unit time interval. The reference scheme's complexity per unit time interval and grid cell is thus $\mathcal{O}\left(\frac{M^2}{\epsilon}\right)$. This computational complexity is problematic for vanishing $\epsilon$.

For a general hierarchical micro-macro acceleration, the runtime is determined by the macro time step $\Delta t$ and the complexity of the four steps of the method: micro step, restriction step, macro step, and matching step. Note that the size of the micro time step $\delta t$ does not influence the runtime, as only a constant number of micro steps with size $\delta t = \mathcal{O}\left(\epsilon\right)$ are performed before each macro step with CFL-size $\Delta t = \mathcal{O}\left(\Delta x\right)$. 
The complexity of the four steps can be approximated as follows:

\begin{itemize}
    \item[1.] The micro step using a moment model of size $M$ has complexity $\mathcal{O}(M^2)$. Additional $K$ micro steps, especially for the PI and CPI methods, increase the complexity to $\mathcal{O}\left((K+1)M^2\right)$.
    \item[2.] The restriction from $M$ to $L$ moments is typically the identity operator and formally has complexity $\mathcal{O}(L)$. 
    \item[3.] The macro step using a macroscopic moment model of size $L$ has complexity $\mathcal{O}(L^2)$. The simple extrapolation of PI and CPI decreases the complexity to $\mathcal{O}(L)$.
    \item[4.] The matching step based on a simple matrix-vector multiplication has complexity $\mathcal{O}\left((M-L)^2\right)$ for non-linear matching like the mMHME or the CPI method. See section \ref{sec:efficient_matching} for a more efficient implementation of the non-linear matching that has complexity $\mathcal{O}(M)$. In case of a linear model or full extrapolation, no matching is necessary and the formal complexity of setting the variables is $\mathcal{O}(M-L)$. 
\end{itemize}

The complexities of the different steps for the mMHME, mMHSM, PI, CPI methods are summarized in table \ref{tab:methods_complexity}. The efficient implementation of the matching step described in section \ref{sec:efficient_matching} reduces the matching runtime  to $\mathcal{O}(M)$. Therefore, the runtime of the restriction step, macro step, and matching step can be neglected with respect to the micro step, which has runtime $\mathcal{O}(M^2)$. Together with a CFL-type time step $\Delta t = \mathcal{O}\left(\Delta x\right)$, a typical mM method has complexity $\mathcal{O}\left(\frac{M^2}{\Delta x}\right)$ per unit time interval and grid cell. In comparison to the reference method with complexity $\mathcal{O}\left(\frac{M^2}{\epsilon}\right)$, the mM method is beneficial for small $\epsilon$ or and not too fine grid size $\Delta x$.

\begin{table}[h]
    \centering
    \begin{tabular}{|c||c|c||c|c|c||c||c|}
       \hline
         & \multicolumn{2}{c||}{micro} & \multicolumn{3}{c||}{macro} & \multicolumn{1}{c||}{match.} & \multicolumn{1}{c|}{restr.} \\ \hline
         & model & $\mathcal{O}$ & model & $L$ & $\mathcal{O}$ & $\mathcal{O}$ &  $\mathcal{O}$ \\ \hline \hline 
        mMHME & HME & $M^2$ & Euler & $3$ & $L^2$ & $(M-L)^2 \, \vee \, M$ & $L$  \\ \hline 
        mMHSM & HSM & $M^2$ & Euler & $3$ & $L^2$ & $M-L$ & $L$  \\ \hline 
        PI & $(K+1)$HME & $(K+1)M^2$ & extrap. & $M$ & $L$ & $M-L$ & $L$  \\ \hline
        CPI & $(K+1)$HME & $(K+1)M^2$ & extrap. & $L \leq M$ & $L$ & $(M-L)^2 \, \vee \, M$ & $L$  \\ \hline
    \end{tabular}
    \caption{Hierarchical micro-macro acceleration choices with computational complexity for a single cell.}
    \label{tab:methods_complexity}
\end{table}

\section{Numerical Results}
\label{sec:Results}

\subsection{Two-beam convergence tests}
\label{sec:2beam}
In a first test case, we investigate the effect of the micro and macro models on the accuracy of the hierarchical micro-macro acceleration with respect to a reference solution. We consider the standard two-beam test case, which was used in \cite{Schaerer2015,Koellermeier2017d,Koellermeier2021} for the HME and HSM models. In \cite{Koellermeier2021} it was shown that a straightforward simulation using the HSM model leads to similar solutions as using the HME model. For this reason we restrict the presentation of results to the HME micro model here. The test case considers the discontinuous initial data
%
\begin{equation}
    \var(0,x) = \left\{
  \begin{array}{ll}
    \var_M^L = \left( 1,0.5,1,0,\ldots,0\right)^T, & \textrm{if } x < 0, \\
    \var_M^R = \left( 1,-0.5,1,0,\ldots,0\right)^T, & \textrm{if } x > 0, \\
  \end{array}
\right.
    \label{e:2beam_IC}
\end{equation}
which models two colliding beams of particles. The test case represents a challenge for moment models as the solution is difficult to represent with a polynomial ansatz. In the collisionless case, an analytical solution exists, which is the sum of two Maxwellians, see \cite{Schaerer2015}. For the simulations, we use a computational grid of $500$ cells on $[-10,10]$, an end time of $t_{\textrm{END}}=0.1$, and a macroscopic time step chosen according to the CFL condition for the macroscopic model resulting in $\Delta t \approx 5 \cdot 10^{-4}$. In the figures that follow, we only show part of the negative domain as the test case is symmetric with respect to $x=0$. The relaxation time $\epsilon = 10^{-4}$ is chosen such that the flow is not too far from equilibrium, but still shows some non-equilibrium behavior, which cannot be modelled by a simple macro model, as identified in \cite{Koellermeier2021}. In addition, it was shown in the same work that this relaxation time yields a time step constraint for a standard explicit method. For smaller relaxation times, non-stiff equilibrium model would suffice, while larger relaxation times can be modelled by non-stiff micro models. The reference solution is computed with an HME moment model using $M=10$ moments from \cite{Koellermeier2021}. Note that all hierarchical micro-macro acceleration method for this test use $K+1=2$ micro iterations per step and employ a third-order FORCE scheme for the spatial discretization \cite{Koellermeier2020}. 

\begin{figure}[htb!]
    \centering
    \begin{subfigures}
    \subfloat[Pressure $p$. \label{fig:2beam_mMHMEmicro_p}
    ]{\includegraphics[width=0.5\linewidth]{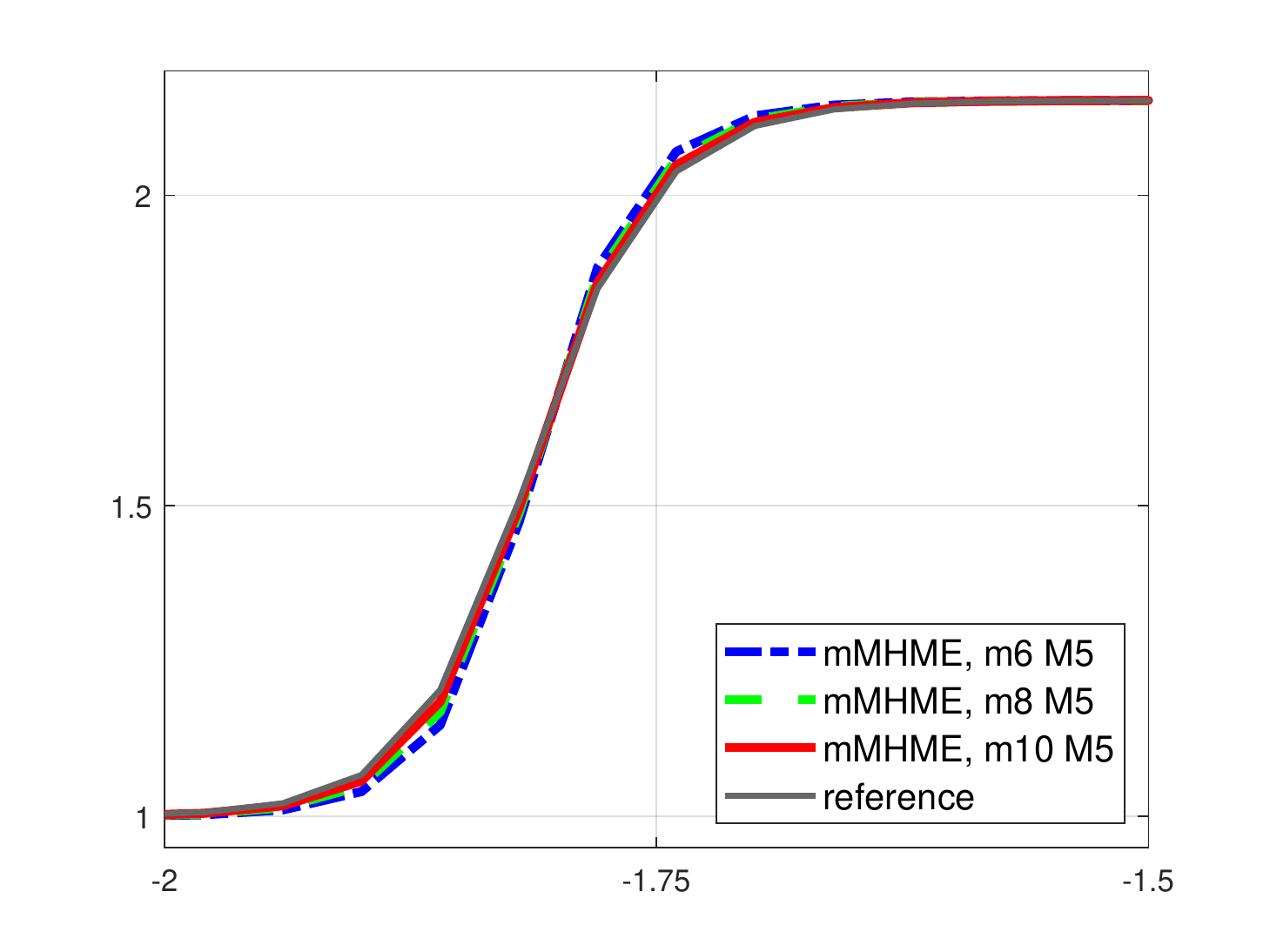}}
    \subfloat[Heat flux $q$. \label{fig:2beam_mMHMEmicro_q}
    ]{\includegraphics[width=0.5\linewidth]{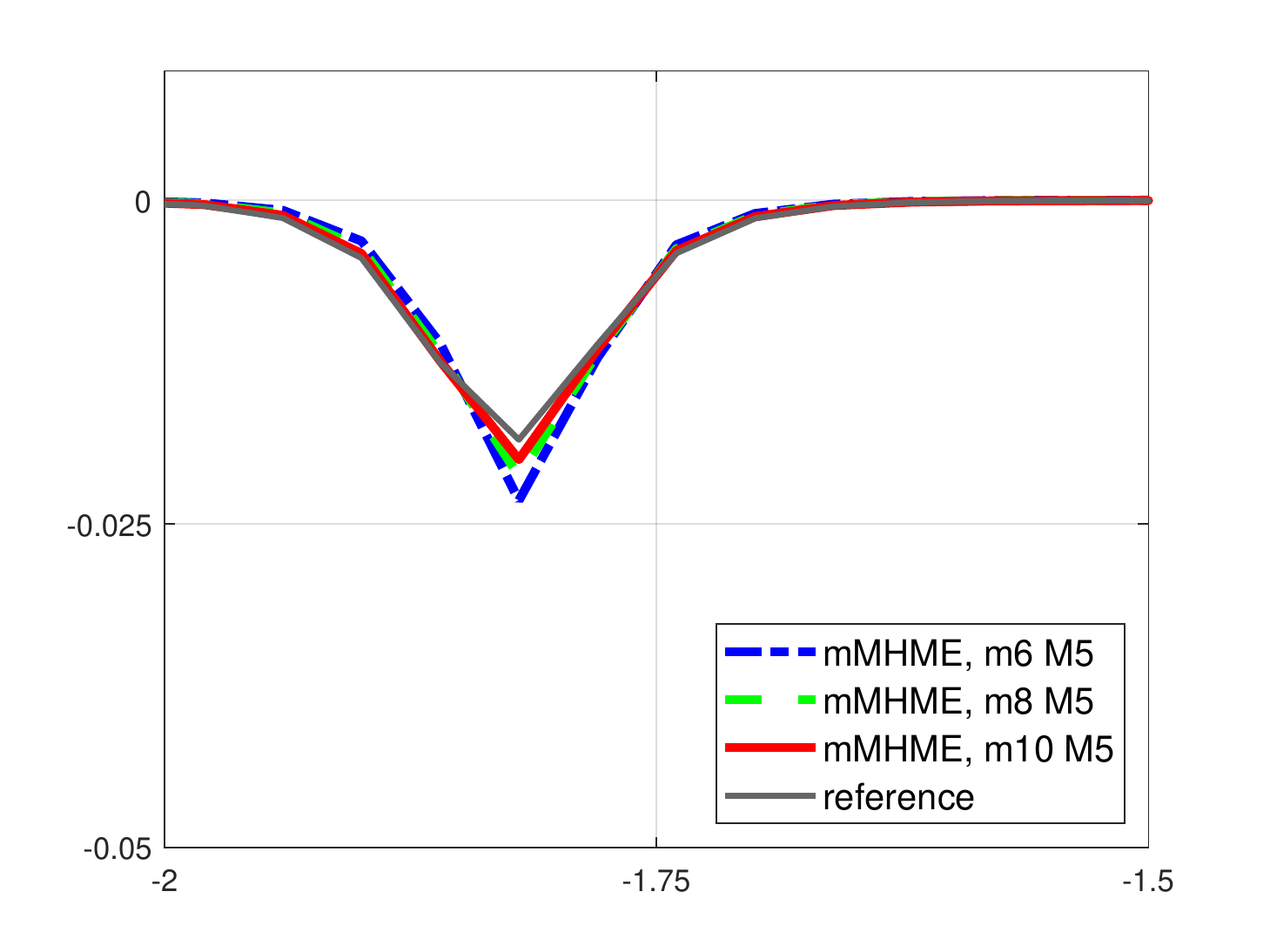}}
    \end{subfigures}
    \caption{Two-beam results with $\epsilon = 10^{-4}$ for mMHME using different micro models $m=6,8,10$ and fixed macro model $M=5$. A more refined micro model increases the accuracy.}
    \label{fig:2beam_mMHMEmicro}
\end{figure}

In Figure \ref{fig:2beam_mMHMEmicro} the mMHME method is used with a fixed $5-$moment macro model and the effect of different micro models is investigated. The results show that a more refined micro model increases the accuracy. This is especially visible for the heat flux, which is a non-vanishing non-equilibrium variable in this setting. Using a sufficiently fine micro model, a macro model using only $M=5$ moments already yields good accuracy with respect to the reference solution.

\begin{figure}[htb!]
    \centering
    \begin{subfigures}
    \subfloat[Pressure $p$. \label{fig:2beam_mMHMEmacro_p}
    ]{\includegraphics[width=0.5\linewidth]{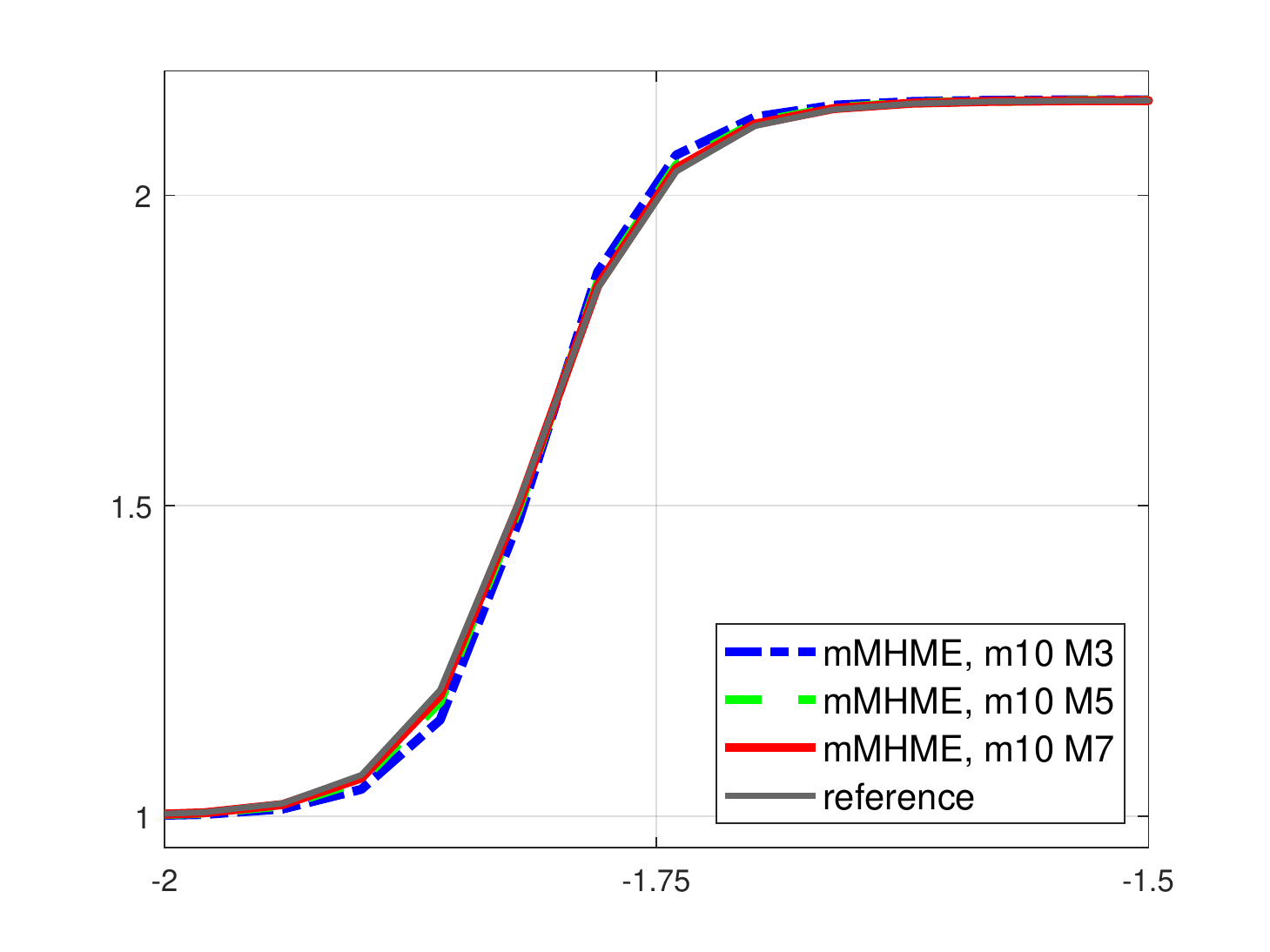}}
    \subfloat[Heat flux $q$. \label{fig:2beam_mMHMEmacro_q}
    ]{\includegraphics[width=0.5\linewidth]{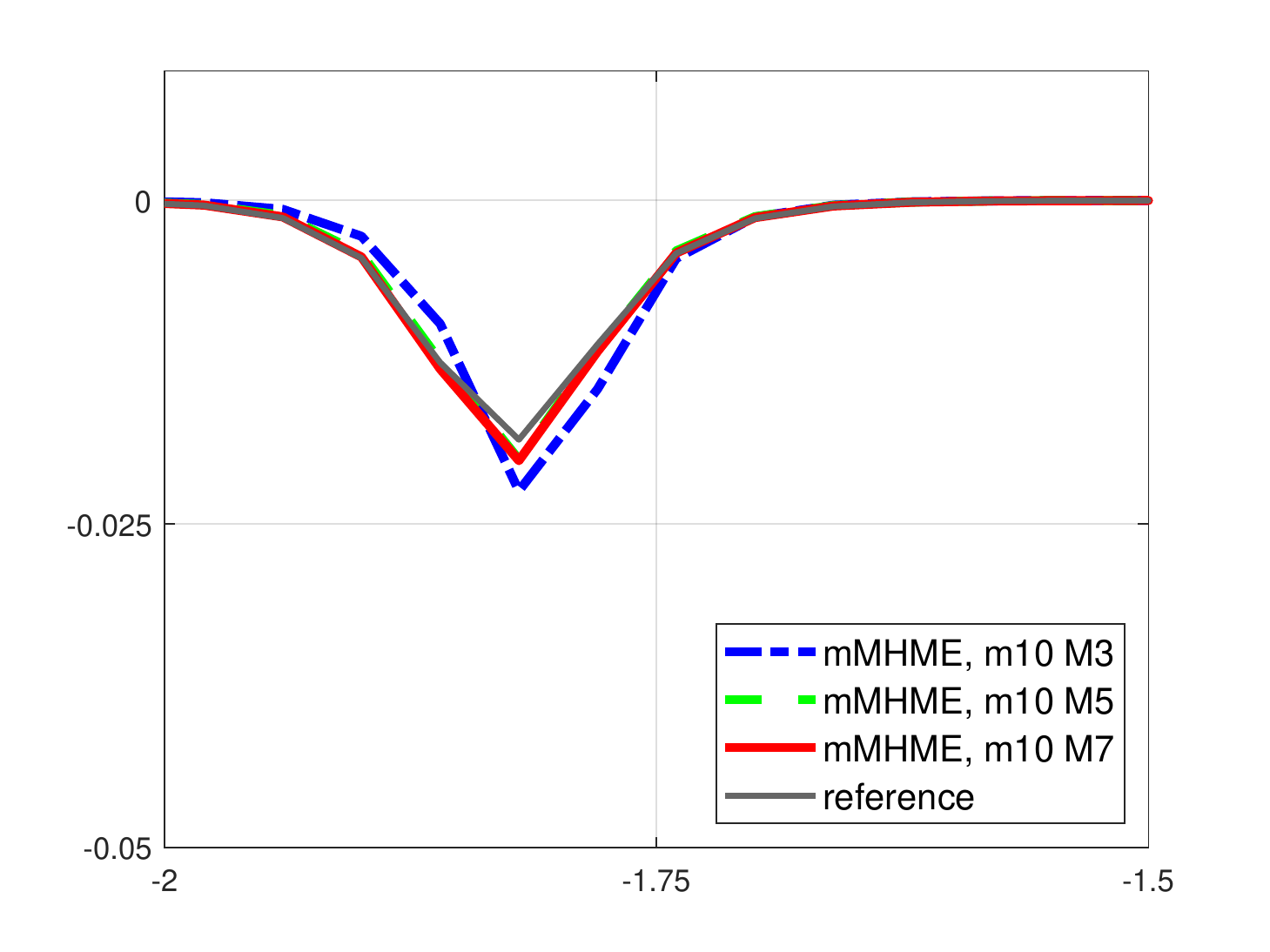}}
    \end{subfigures}
    \caption{Two-beam results with $\epsilon = 10^{-4}$ for mMHME using different macro models $M=3,5,7$ and fixed micro model $m=10$. A more refined macro model increases the accuracy.}
    \label{fig:2beam_mMHMEmacro}
\end{figure}

In Figure \ref{fig:2beam_mMHMEmacro} the mMHME method is used with a fixed $10-$moment micro model and the effect of using different macro models is investigated. Similarly to the micro model, also a more refined macro model increases the accuracy of the simulation. Due to the accurate micro model, the solution quality is already very good, as can be seen for the pressure $p$ in \ref{fig:2beam_mMHMEmacro_p}.

As a second example for the methods of this paper, we consider the CPI method for the same settings as presented above. In Figure \ref{fig:2beam_mMHMEmicro}, the micro model is refined for a constant macro model. Note that the macro model here consists of an extrapolation over a large time step and the number of extrapolated variables is chosen as $L=5$ here. The results show that the the CPI method benefits even more from refining the micro model. The accuracy is slightly worse than for the mMHME model, which can be explained by the more accurate macro step using a fully non-linear model in case of mMHME.
\begin{figure}[htb!]
    \centering
    \begin{subfigures}
    \subfloat[Pressure $p$. \label{fig:2beam_mMCPImicro_p}
    ]{\includegraphics[width=0.5\linewidth]{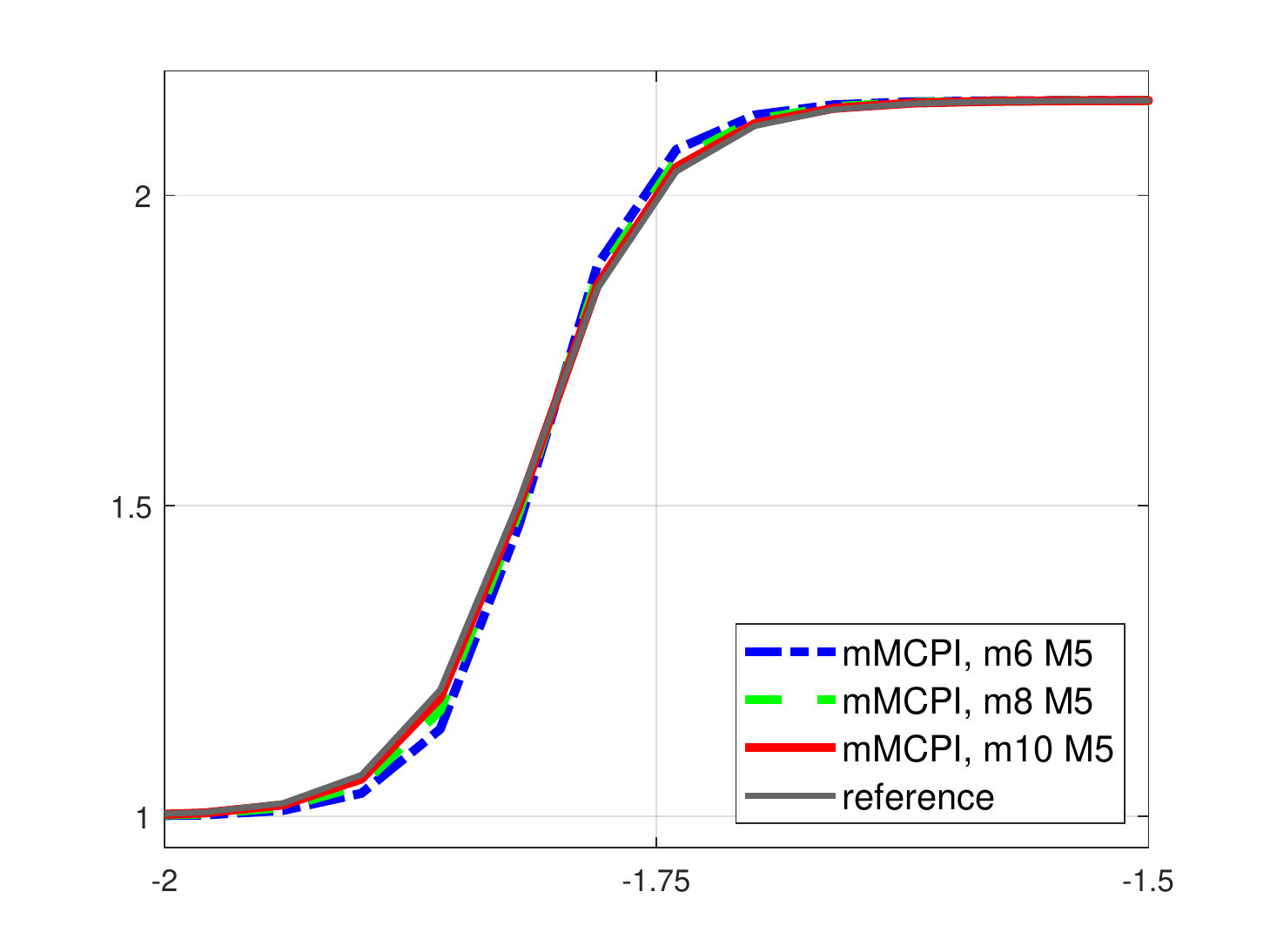}}
    \subfloat[Heat flux $q$. \label{fig:2beam_mMCPImicro_q}
    ]{\includegraphics[width=0.5\linewidth]{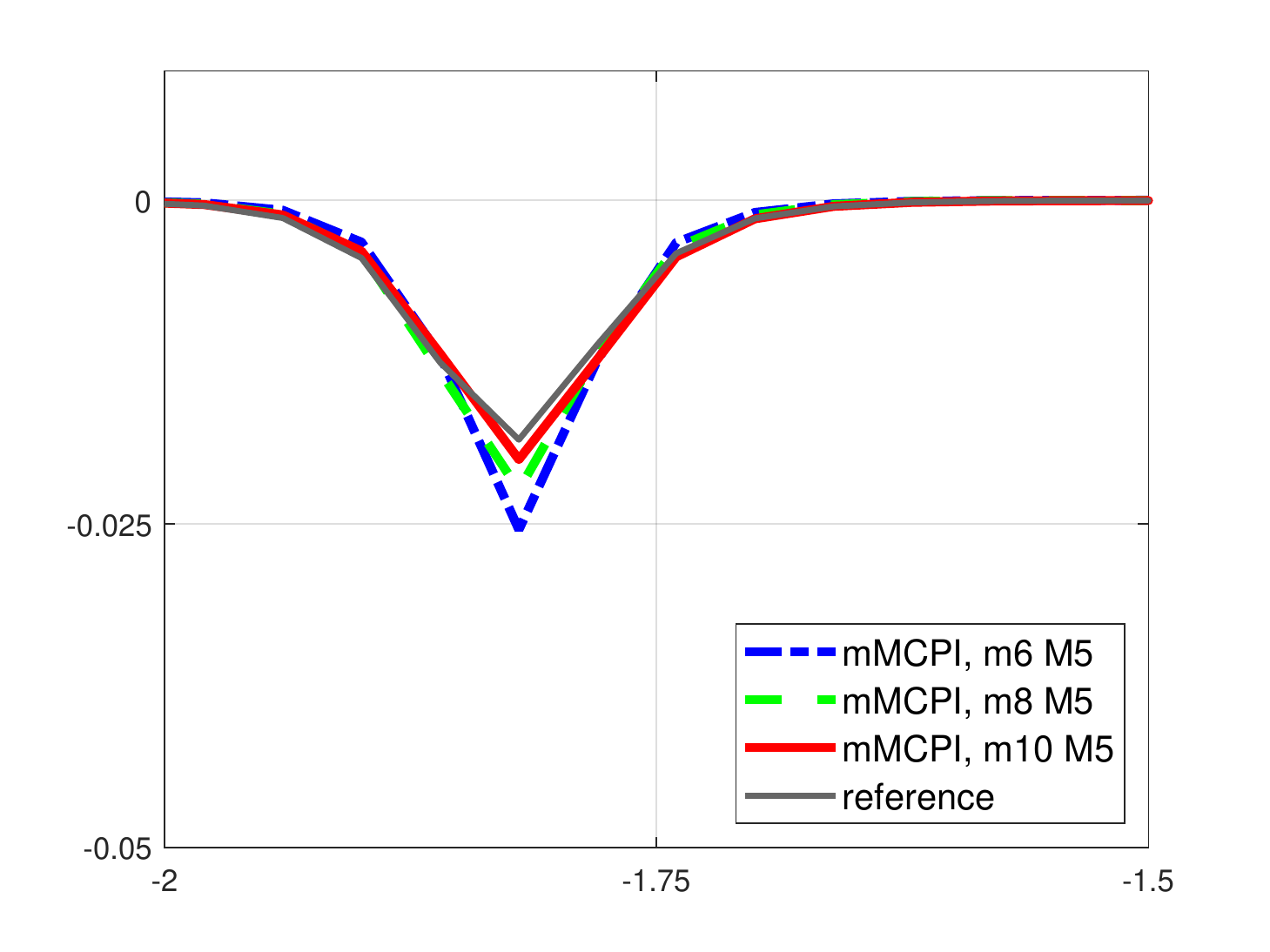}}
    \end{subfigures}
    \caption{Two-beam results with $\epsilon = 10^{-4}$ for CPI using different micro models $m=6,8,10$. A more refined micro model increases the accuracy of Coarse Projective Integration.}
    \label{fig:2beam_mMCPImicro}
\end{figure}

\begin{figure}[htb!]
    \centering
    \begin{subfigures}
    \subfloat[Pressure $p$. \label{fig:2beam_mMCPImacro_p}
    ]{\includegraphics[width=0.5\linewidth]{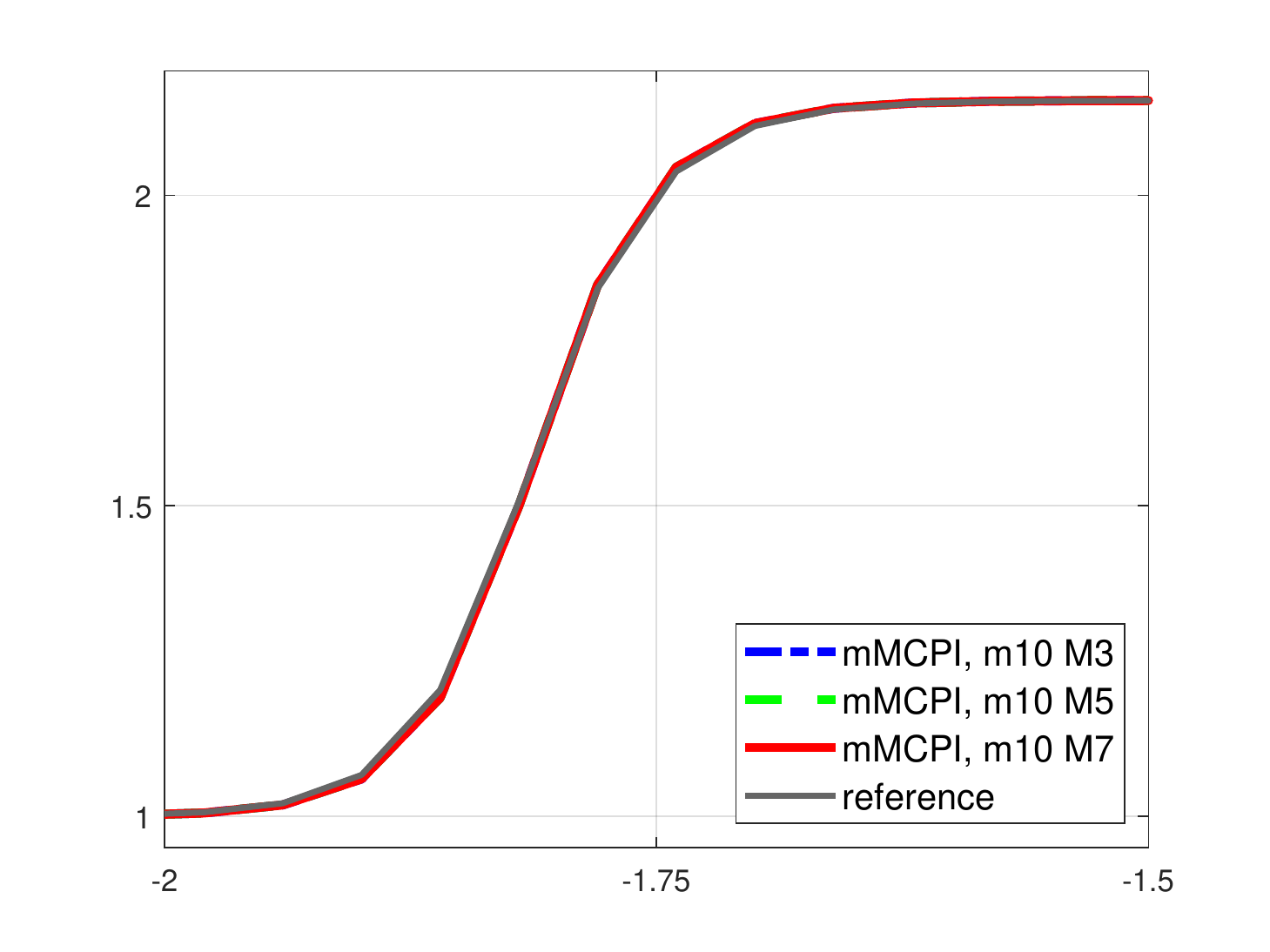}}
    \subfloat[Heat flux $q$. \label{fig:2beam_mMCPImacro_q}
    ]{\includegraphics[width=0.5\linewidth]{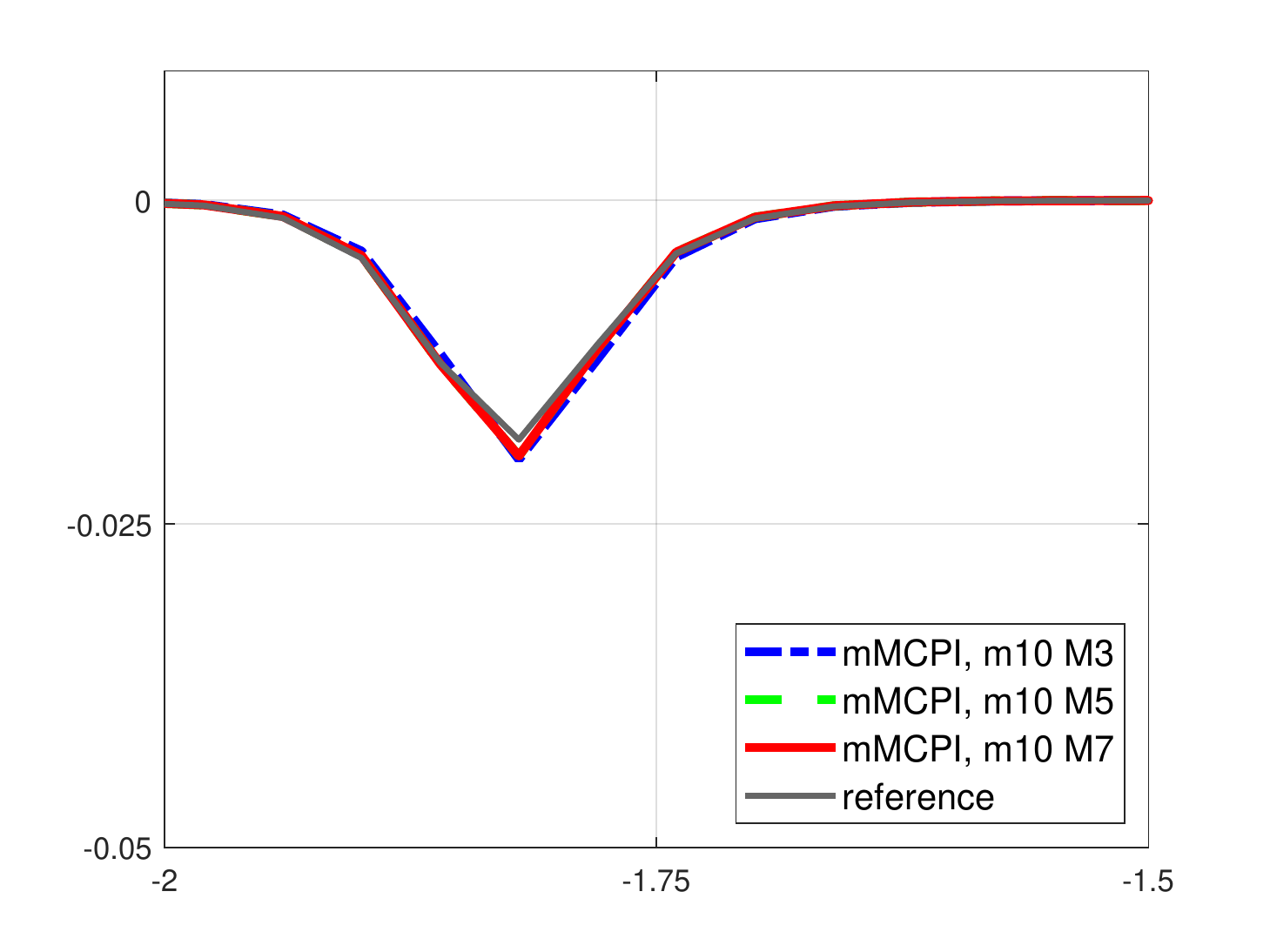}}
    \end{subfigures}
    \caption{Two-beam results with $\epsilon = 10^{-4}$ for CPI using different macro models $M=3,5,7$. Using more macro variables for the Coarse Projective Integration does not significantly increase the accuracy.}
    \label{fig:2beam_mMCPImacro}
\end{figure}

Figure \ref{fig:2beam_mMCPImacro} shows the effect of using different number of moments for the extrapolation step of the CPI method. In contrast to the mMHME method, the different ways of executing the macro step do not influence the accuracy of the CPI method as much. Already using a coarse macro model of only $L=3$ variables leads to a very accurate solution. This indicates that the choice of the micro model is indeed more relevant than the macro model and a coarse macro model might suffice once the micro model has damped the fast modes. 

\begin{figure}[htb!]
    \centering
    \begin{subfigures}
    \subfloat[mMHME, pressure $p$. 
    ]{\includegraphics[width=0.5\linewidth]{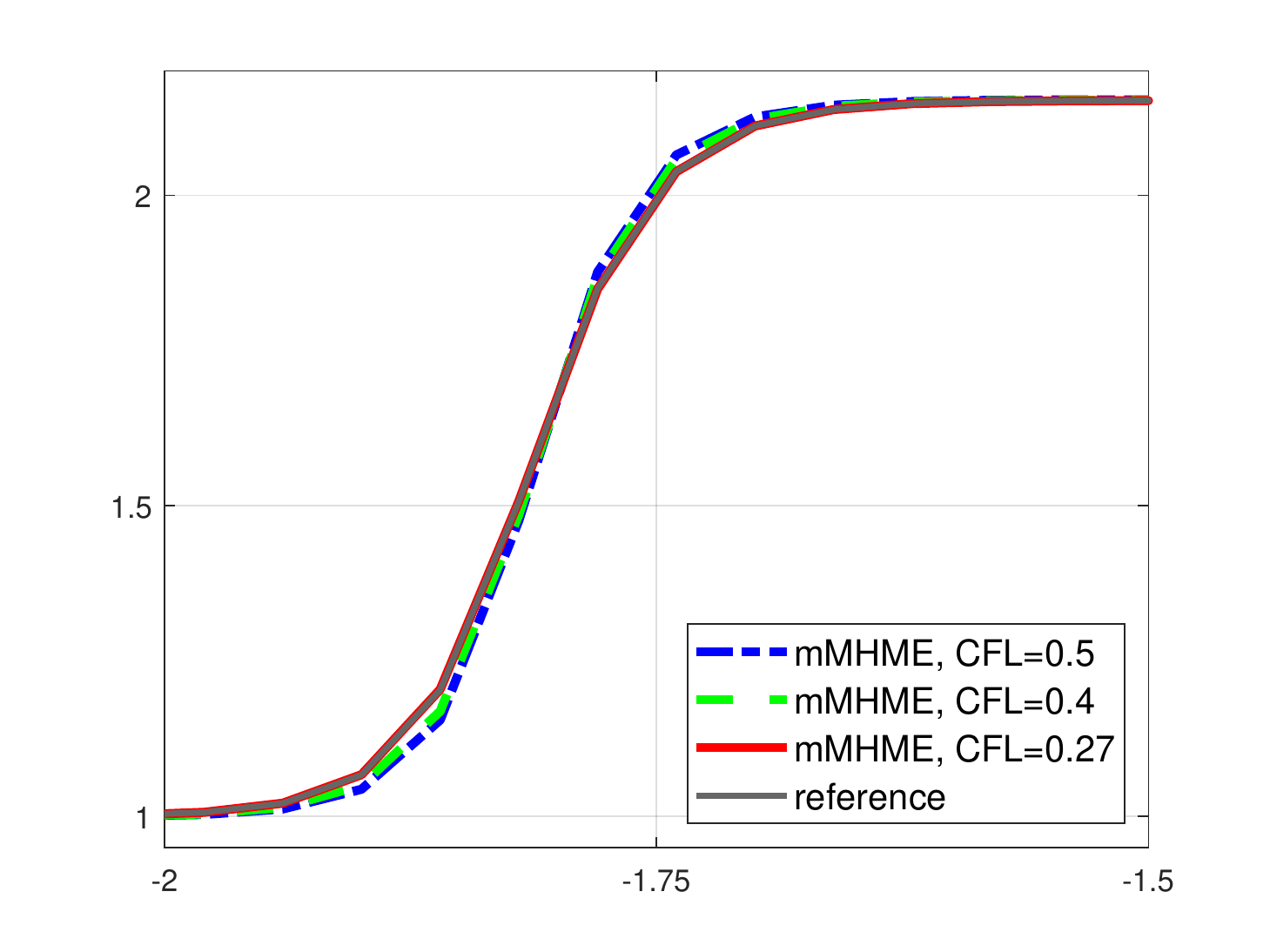}}
    \subfloat[mMHME, heat flux $q$. 
    ]{\includegraphics[width=0.5\linewidth]{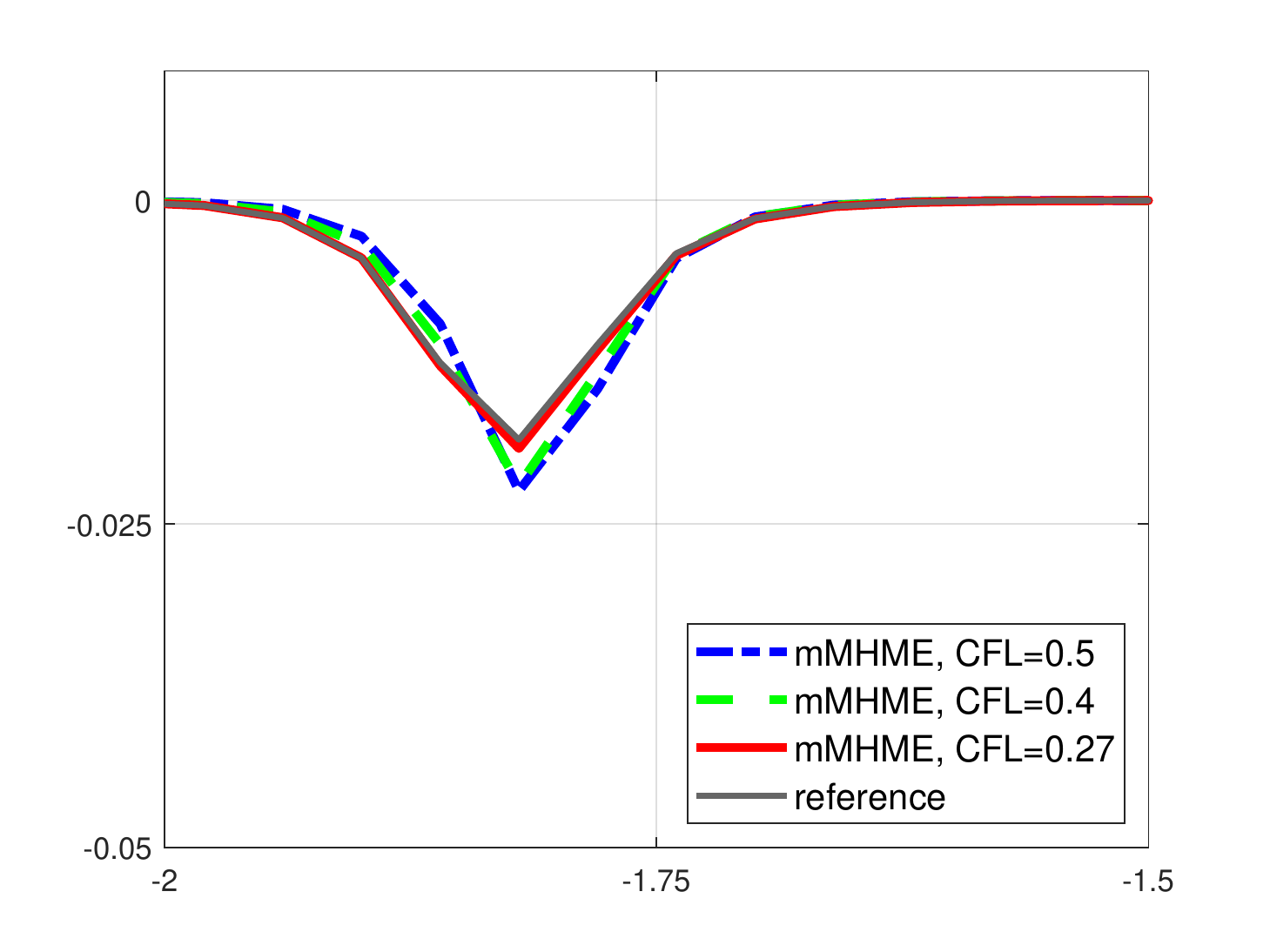}} \\
    \subfloat[CPI, pressure $p$. 
    ]{\includegraphics[width=0.5\linewidth]{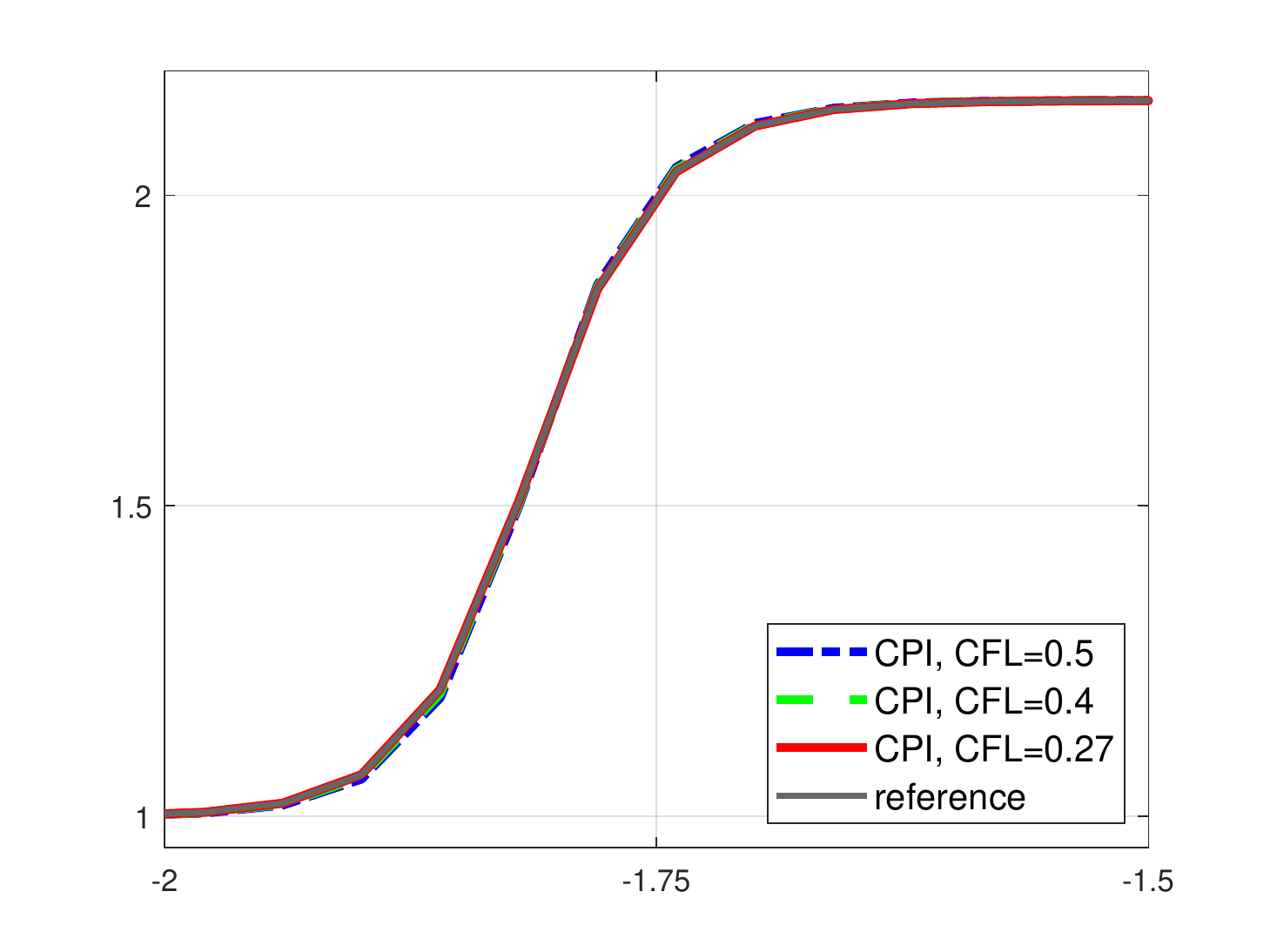}}
    \subfloat[CPI, heat flux $q$.
    ]{\includegraphics[width=0.5\linewidth]{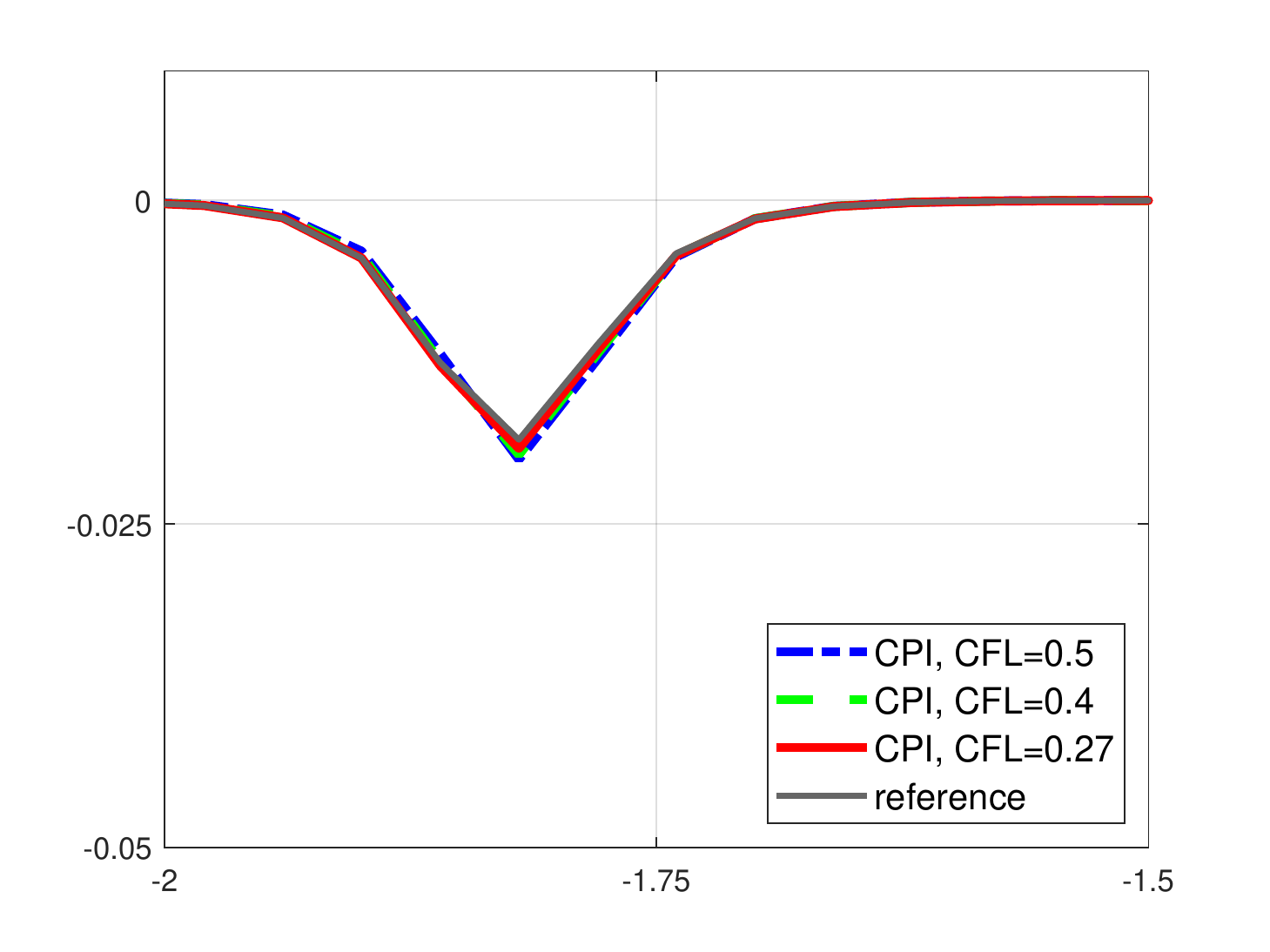}}
    \end{subfigures}
    \caption{Two-beam results with $\epsilon = 10^{-4}$ for mMHME (top row) and CPI (bottom row) using fixed macro model $M=3$ and micro model $m=10$ and different $\Delta t$ given by $CFL=0.5, 0.4, 0.27$ and third order FORCE method. The mM solution converges to the micro solution.}
    \label{fig:2beam_consistency}
\end{figure}

A numerical check of the consistency property, is shown in the bottom row of Figure \ref{fig:2beam_consistency} for the mMHME method. The macro step size is varied by means of a different CFL number $\frac{\Delta t}{\Delta x}$, so that the remaining macro time interval after the micro iterations is decreasing. Both the pressure and heat flux converge to the reference solution, which can be seen as the micro solution. This is numerical evidence for the consistency property, which was analytically derived in section \ref{sec:Consistency_proof}.

The same consistency property also holds numerically for the CPI method, which uses extrapolation of the first $M$ moments. The results are presented in the bottom row of Figure \ref{fig:2beam_consistency}. A decrease of the macro time step size by means of a smaller CFL number also leads to convergence towards the reference micro solution. This shows that the hierarchical micro-macro acceleration is consistent with the underlying micro model.

\subsection{Runtime comparison}
\label{sec:speedup}
For the runtime comparison of the 1D case, we consider the setup of the previous two-beam test case and the following three methods for comparison:
\begin{itemize}
    \item[1.] a standard macro model with $L=3$ equations (Euler equations) exhibiting no stiffness constraints and using a CFL-type time step of $\Delta t = 5 \cdot 10^{-4}$.
    \item[2.] a standard micro model with $M=10$ equations (moment model) exhibiting stiffness constraints such that $\Delta t \approx \epsilon$ for small $\epsilon$. 
    \item[3.] the new hierarchical micro-macro acceleration from section \ref{sec:3_mMHME} with $M=10$ moments for the micro model and $L=3$ equations for the macro model (Euler equations). We perform $K+1 = 2$ steps of the micro model. The computational overhead occurring from the matching, restriction, and macro step is neglected as discussed in section \ref{sec:Complexity}. The method overcomes the stiffness constraints and uses a CFL-type time step of $\Delta t = 5 \cdot 10^{-4}$.
\end{itemize}

We expect that the macro model results in a fast but inaccurate solution, while the full micro model results in a very slow but accurate solution. The hierarchical micro-macro acceleration should be much faster than a full micro model while retaining the accuracy close to equilibrium as shown in section \ref{sec:2beam}.

In table \ref{tab:methods_runtime} we compare the speedup as computed using the complexity analysis in section \ref{sec:Complexity} and take the micro model as reference. (Note that we use the complexity analysis as actual runtime largely depends on the implementation, especially for a fast 1D code. For a comparison of the actual runtimes for the 2D test case, we refer to Figure \ref{tab:methods_runtime_FFS}.) A hierarchical micro-macro acceleration only makes sense to use for small values of $\epsilon$ that result in a time step size restriction due to stiffness of the micro model, here from $\epsilon=10^{-4}$ on. We know that the runtime of a full micro simulation is increasing quickly with smaller relaxation time $\epsilon$, due to the more severe time step constraint, while the runtime of the hierarchical micro-macro acceleration in this case does not depend on $\epsilon$. The hierarchical micro-macro acceleration thus yields a significant speedup of up to $255.5$ in comparison to the full micro solution while still resulting in accurate solutions as shown in the previous test case. We remind the reader that the macro model (Euler equations) results in a large speedup, but lacks accuracy.
We note that for more complex test cases, e.g., non-linear collision operators, the number of micro steps may need to be increased to satisfy the stability conditions of hierarchical micro-macro acceleration, see \cite{Koellermeier2021}. 

\begin{table}[h]
    \centering
    \begin{tabular}{|c||c|c|c|c|c|}
       \hline
        method & complexity & $\epsilon = 10^{-3}$ & $\epsilon = 10^{-4}$ & $\epsilon = 10^{-5}$ & $\epsilon = 10^{-6}$ \\ \hline 
       \hline
        macro & $\mathcal{O}\left(\frac{L^2}{\Delta t}\right)$ & 11.1 & 55.5 & 555.5 & 5555.5  \\ \hline 
        micro-macro & $\mathcal{O}\left((K+1)\frac{M^2}{\Delta t}\right)$ & - & 2.5 & 25.5 & 255.5  \\ \hline 
        micro & $\mathcal{O}\left(\frac{M^2}{\epsilon}\right)$ & 1 & 1 & 1 & 1  \\ \hline
    \end{tabular}
    \caption{Two-beam test case speedup with respect to micro computation for comparison between macro and micro-macro methods for $L=3$, $M=10$, $K=1$, CFL-type $\Delta t=5\cdot 10^{-4}$ and varying values for $\epsilon$.}
    \label{tab:methods_runtime}
\end{table}


\subsection{2D forward facing step application}
\label{sec:2D}
The 2D application case is a rarefied supersonic flow over a forward facing step. It has been studied, among others in \cite{Bogolepov1983,Stueer1999} and for the hyperbolic moment models in \cite{Koellermeier2017d,Koellermeier2021}. At the inlet of a rectangular domain the flow enters with a Mach number $\textrm{Ma} = 3$ and a step close to the inlet generates shock waves, and subsequent separation, reattachment, and reflection of the shock wave, see \cite{Koellermeier2017d} for details.

\begin{figure}[htbp!]
    \centering
    \begin{subfigures}
    \subfloat[Euler equations, $\Delta t = 0.001$. \label{fig:0p01EulerKn0p01deltaT0p001t11scale}
    ]{\begin{overpic}[width=0.95\textwidth]{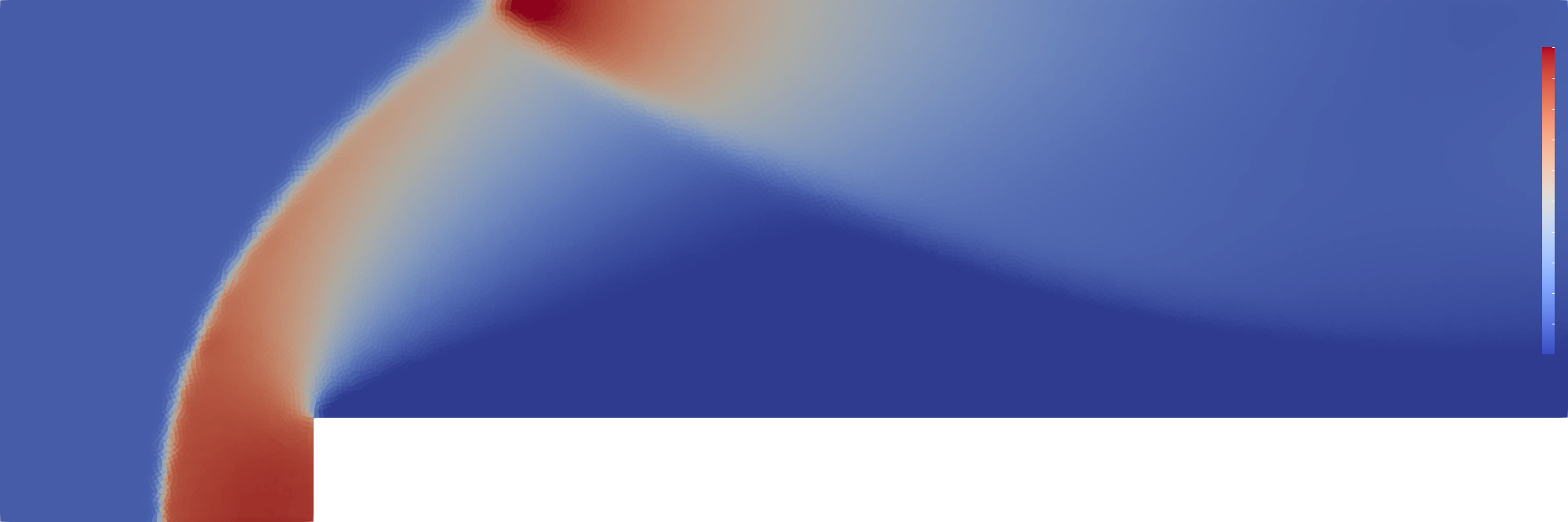}
        \put(94,29.3){$6.5$}
        \put(94,10){$1.5$}
        \put(10.3,0){\color{black}\line(0,1){33.3}}
        \put(31.7,33.3){\color{black}\line(0,-1){26.6}}
    \end{overpic}}\\
    \subfloat[Accurate reference model: HME, $M=3$, splitting scheme $\Delta t = 0.001$. \label{fig:0p01HMEM3Kn0p001deltaT0p001t11}
    ]{\begin{overpic}[width=0.95\textwidth]{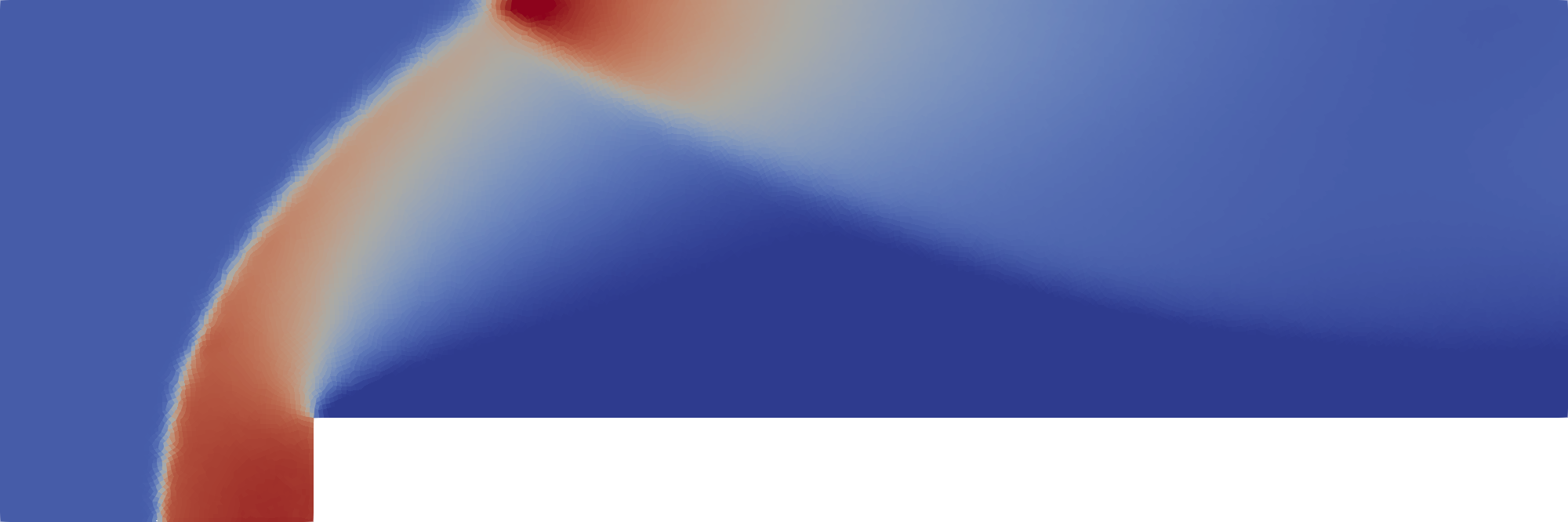}
        \put(10.3,0){\color{black}\line(0,1){33.3}}
        \put(31.7,33.3){\color{black}\line(0,-1){26.6}}
    \end{overpic}}\\
    \subfloat[Smeared out micro model, HME, $M=3$, $\Delta t = 0.0001$. \label{fig:0p01HMEM3Kn0p001deltaT0p0001t11}
    ]{\begin{overpic}[width=0.95\textwidth]{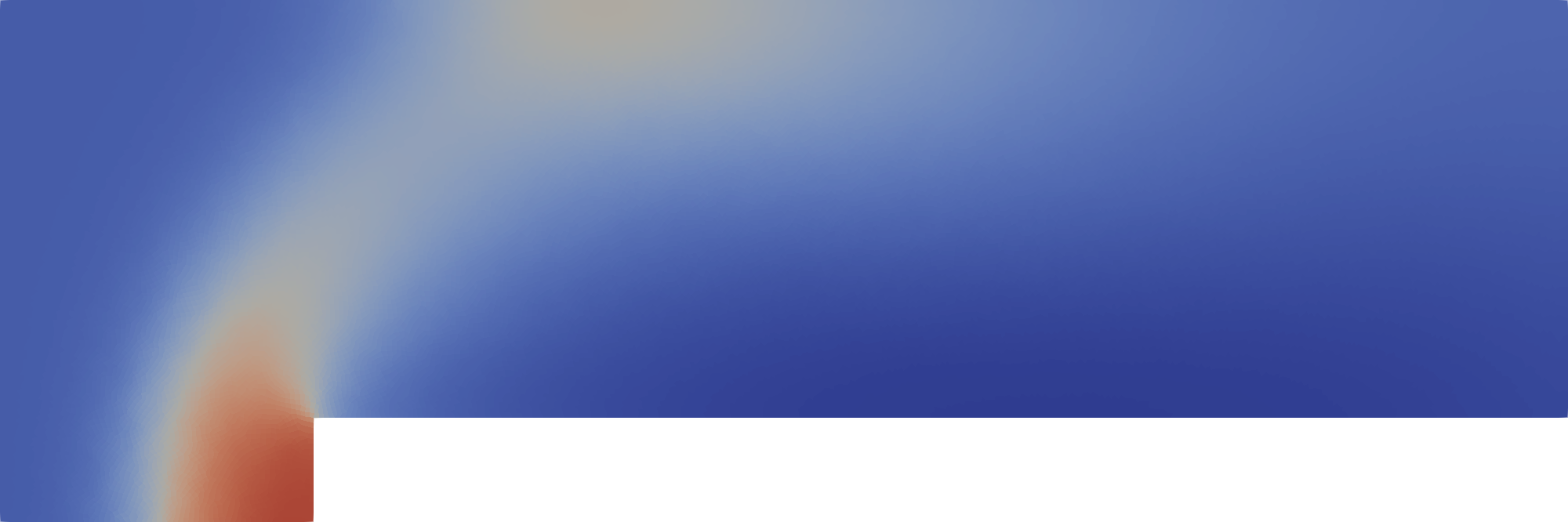}
        \put(10.3,0){\color{black}\line(0,1){33.3}}
        \put(31.7,33.3){\color{black}\line(0,-1){26.6}}
    \end{overpic}}
    \end{subfigures}
    \caption{Forward facing step horizontal velocity $\vel_x$ for $\epsilon = 10^{-3}$, $t_{\textrm{end}}=1$ for reference Euler (a), micro splitting method (b), and diffusive micro model (c). Black lines marking the bow shock and reflection point of the reference micro solution for comparison.}
    \label{fig:ffs_ref}
\end{figure}
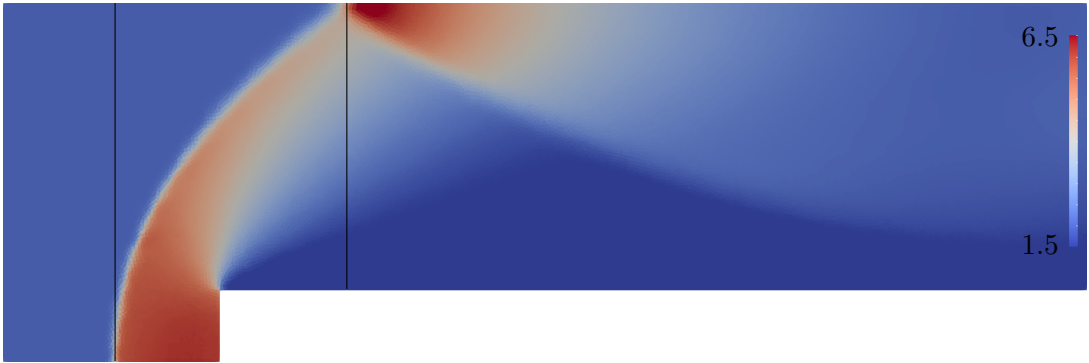
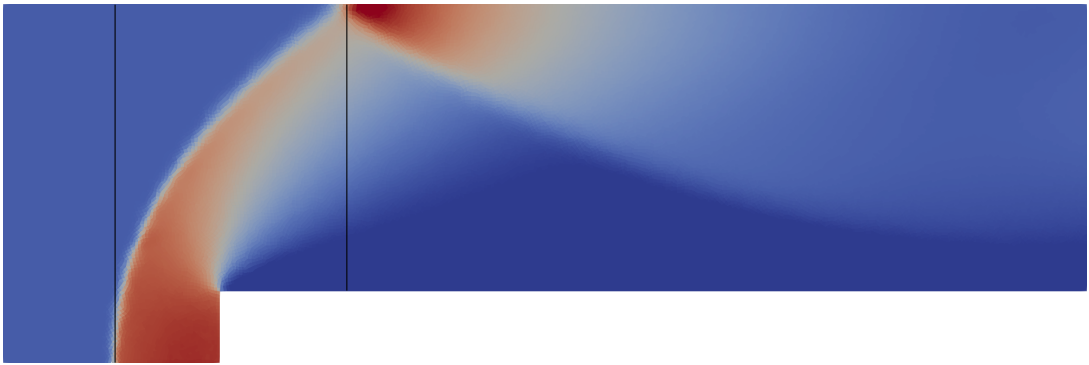
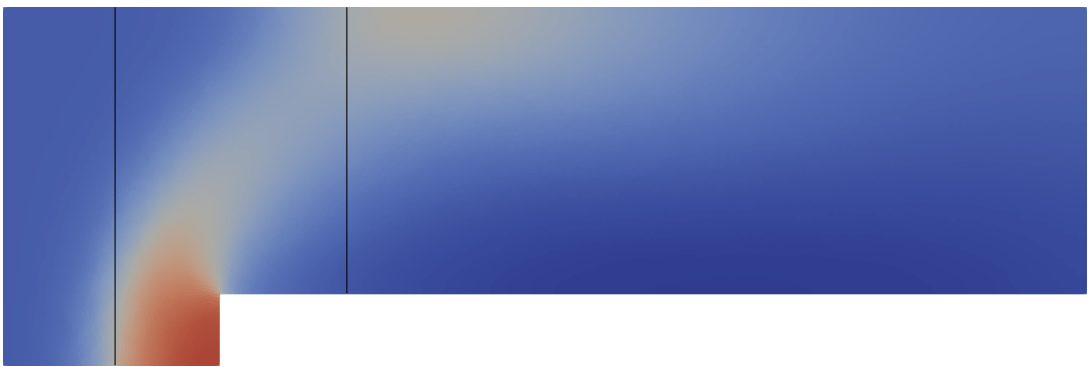

\begin{figure}[htbp!]
    \centering
    \begin{subfigures}
    \subfloat[Projective Integration, HME with $M=3$, $\delta t = 10^{-4}$, $K=3$, $\Delta t = 0.001$. \label{fig:0p01HMEM3Kn0p001deltaT0p001InnerdeltaT0p0001K3Extrapt11}
    ]{\begin{overpic}[width=0.95\textwidth]{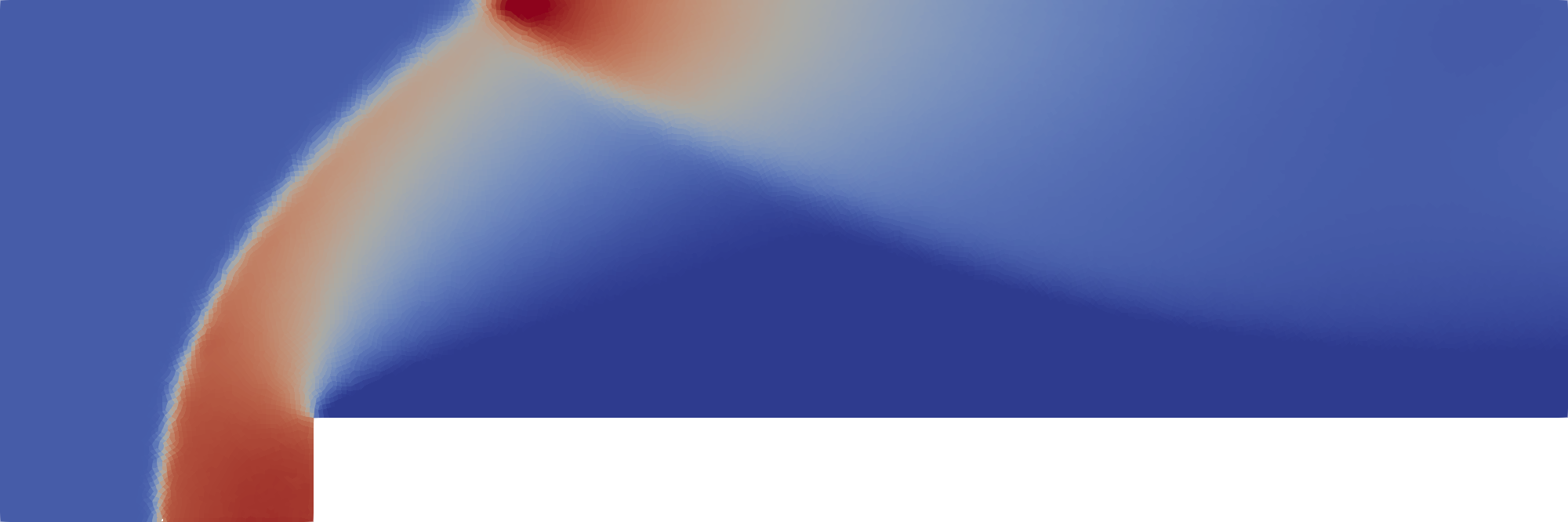}
        \put(10.3,0){\color{black}\line(0,1){33.3}}
        \put(31.7,33.3){\color{black}\line(0,-1){26.6}}
    \end{overpic}}\\
    \subfloat[Coarse Projective Integration, HME with $M=3$, $\delta t = 10^{-4}$, $K=1$, $\Delta t = 0.001$. \label{fig:0p01HMEM3Kn0p001deltaT0p001InnerdeltaT0p0001K1Coarset11}
    ]{\begin{overpic}[width=0.95\textwidth]{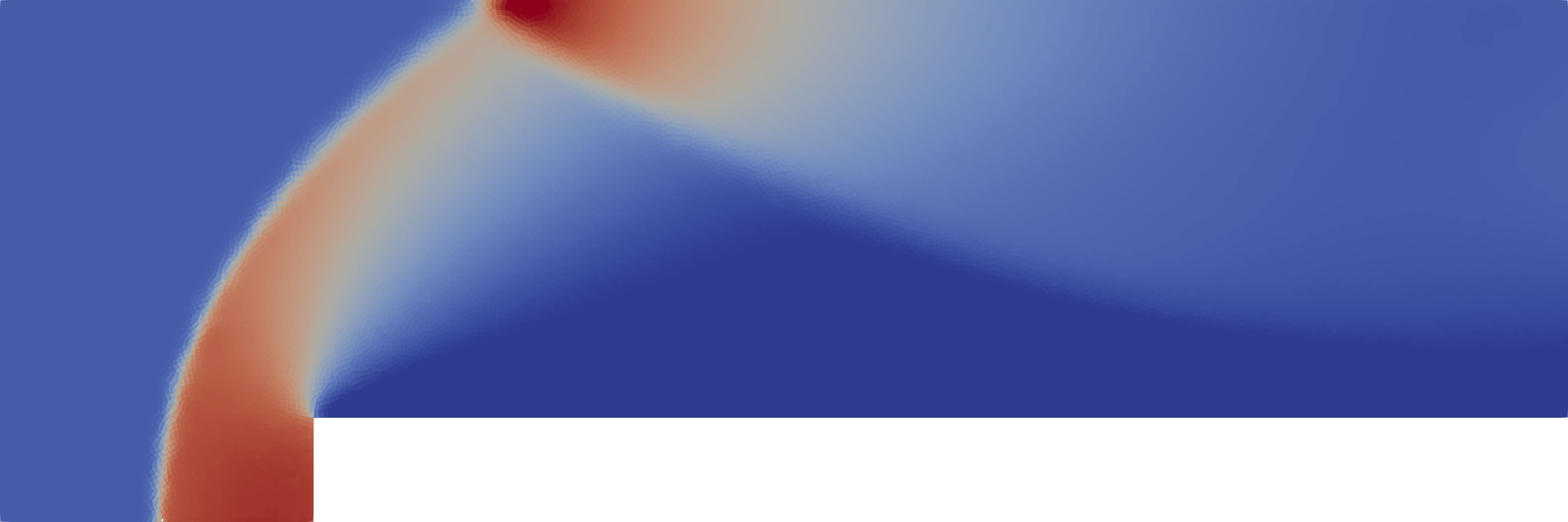}
        \put(10.3,0){\color{black}\line(0,1){33.3}}
        \put(31.7,33.3){\color{black}\line(0,-1){26.6}}
    \end{overpic}}\\
    \subfloat[Micro-macro method, HME with $M=3$, $\delta t = 10^{-4}$, $K=1$, macro Euler with $\Delta t = 0.001$. \label{fig:0p01HMEM3Kn0p001deltaT0p001InnerdeltaT0p0001K1t11}
    ]{\begin{overpic}[width=0.95\textwidth]{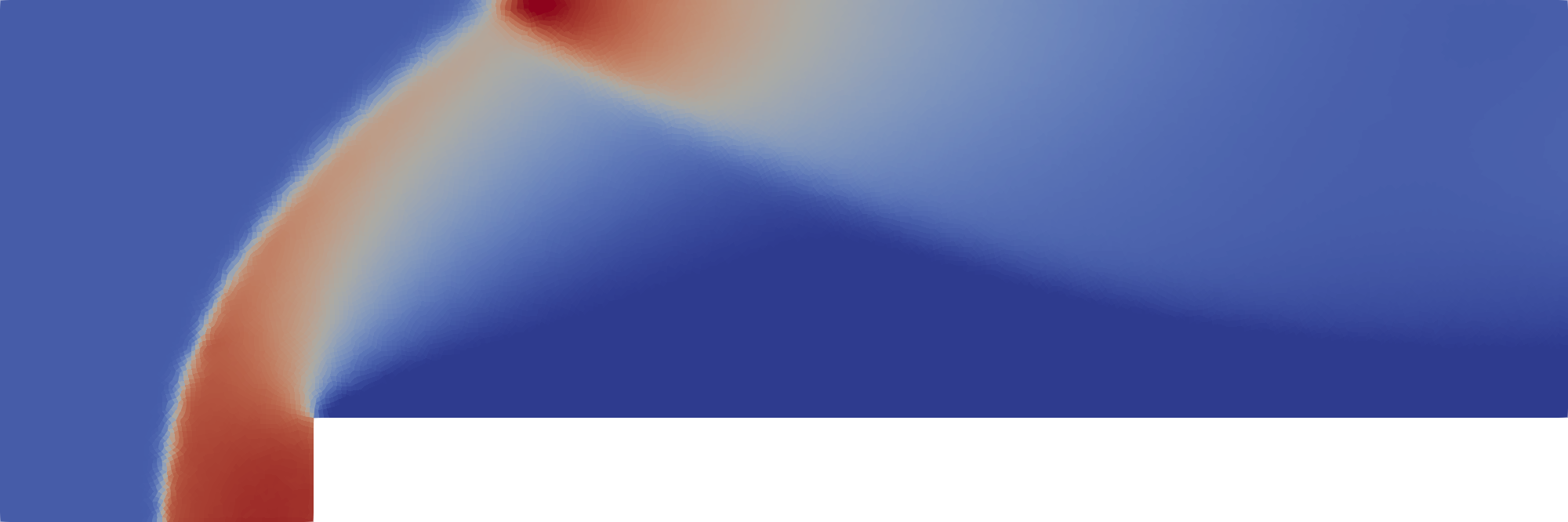}
        \put(10.3,0){\color{black}\line(0,1){33.3}}
        \put(31.7,33.3){\color{black}\line(0,-1){26.6}}
    \end{overpic}}
    \end{subfigures}
    \caption{Forward facing step horizontal velocity $\vel_x$ for $\epsilon = 10^{-3}$, $t_{\textrm{end}}=1$ for Projective Integration (a), Coarse Projective Integration (b), and micro-macro method (c). Black lines marking the bow shock and reflection point of the reference micro solution from Figure \ref{fig:ffs_ref} for comparison.}
    \label{fig:ffs}
\end{figure}
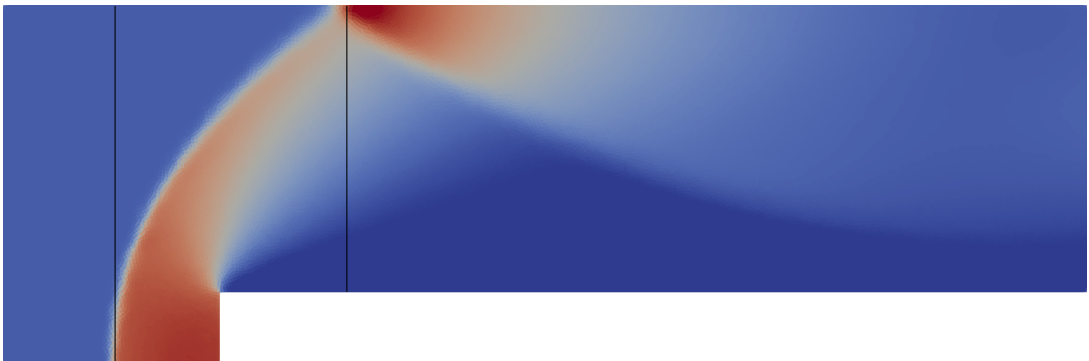
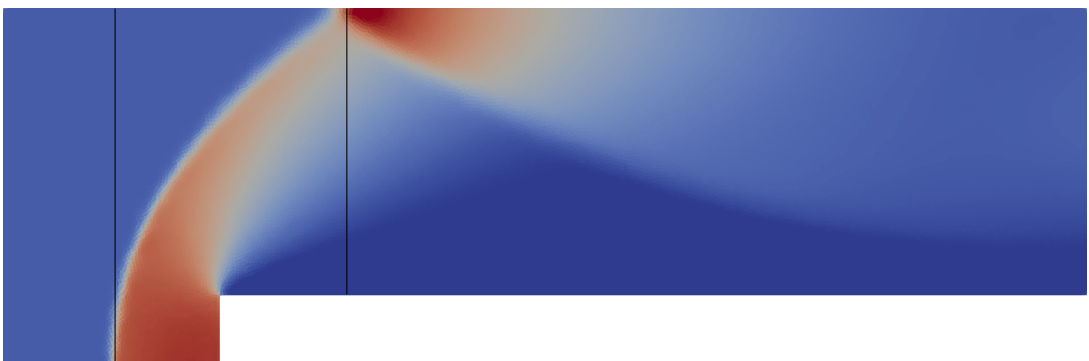
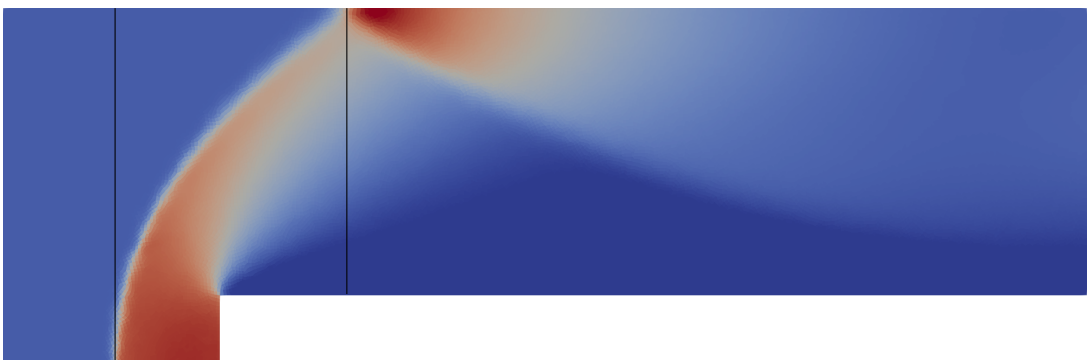

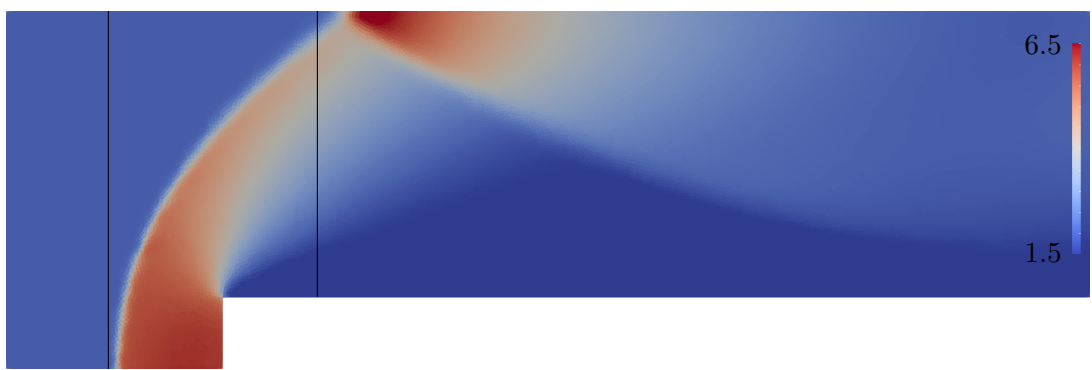
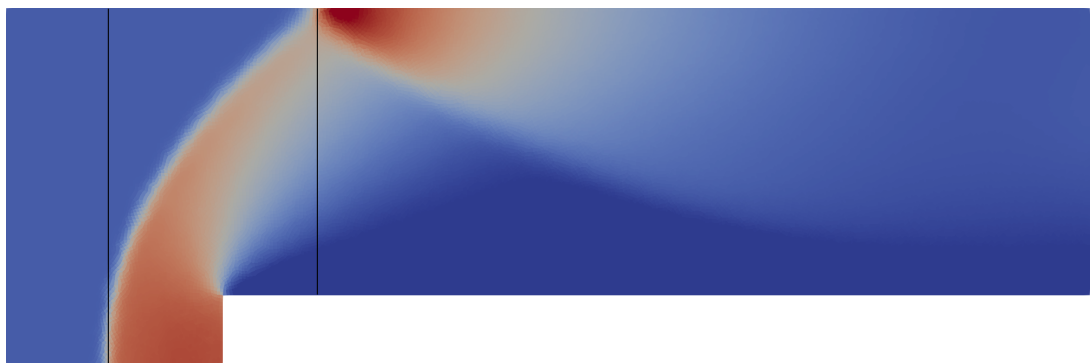
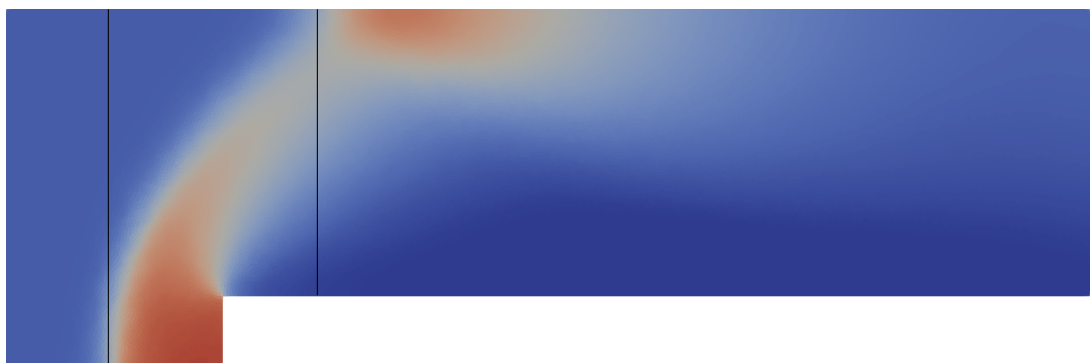
\begin{figure}[htbp!]
    \centering
    \begin{subfigures}
    \subfloat[Euler equations, $\Delta t = 0.001$. \label{fig:0p01EulerKn0p01deltaT0p001t11scale2}
    ]{\begin{overpic}[width=0.95\textwidth]{figures/0p01EulerKn0p01deltaT0p001t11scale}
        \put(94,29.3){$6.5$}
        \put(94,10){$1.5$}
        \put(9.4,0){\color{black}\line(0,1){33.3}}
        \put(28.7,33.3){\color{black}\line(0,-1){26.6}}
    \end{overpic}}\\
    \subfloat[Reference model: HME, $M=3$, splitting scheme $\Delta t = 0.001$. \label{fig:0p01HMEM3Kn0p01deltaT0p001splitt11}
    ]{\begin{overpic}[width=0.95\textwidth]{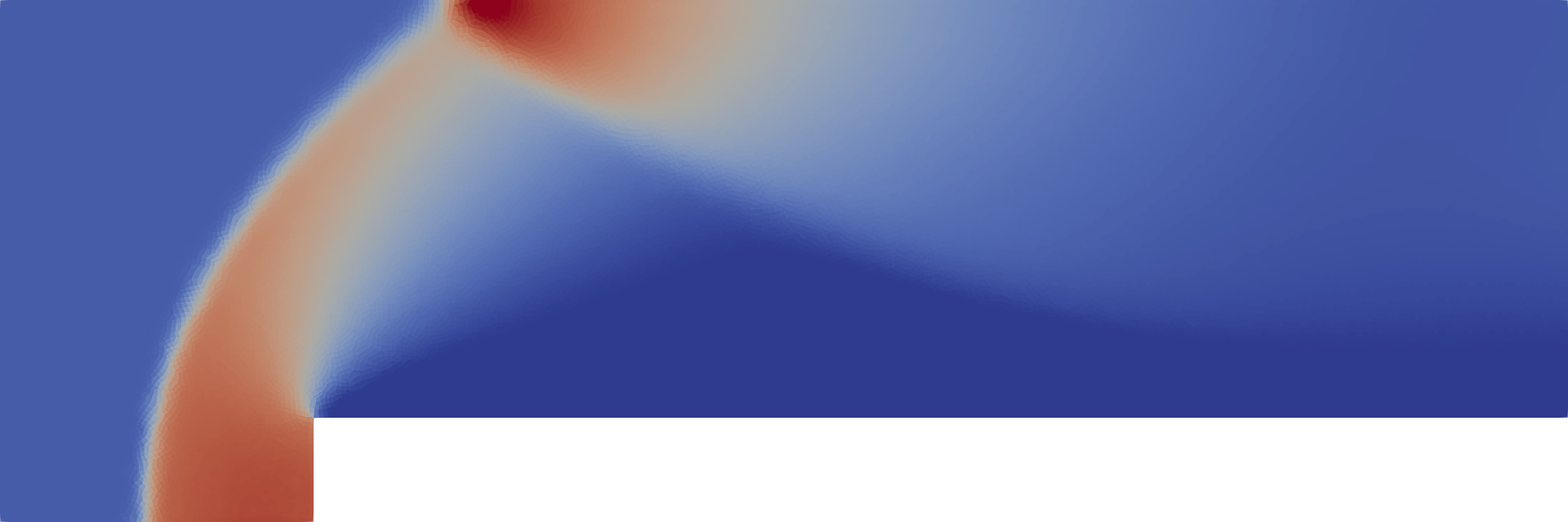}
        \put(9.4,0){\color{black}\line(0,1){33.3}}
        \put(28.7,33.3){\color{black}\line(0,-1){26.6}}
    \end{overpic}}\\
    \subfloat[Micro model, HME, $M=3$, $\Delta t = 2.5\cdot 10^{-4}$. \label{fig:0p01HMEM3Kn0p01deltaT0p00025t11}
    ]{\begin{overpic}[width=0.95\textwidth]{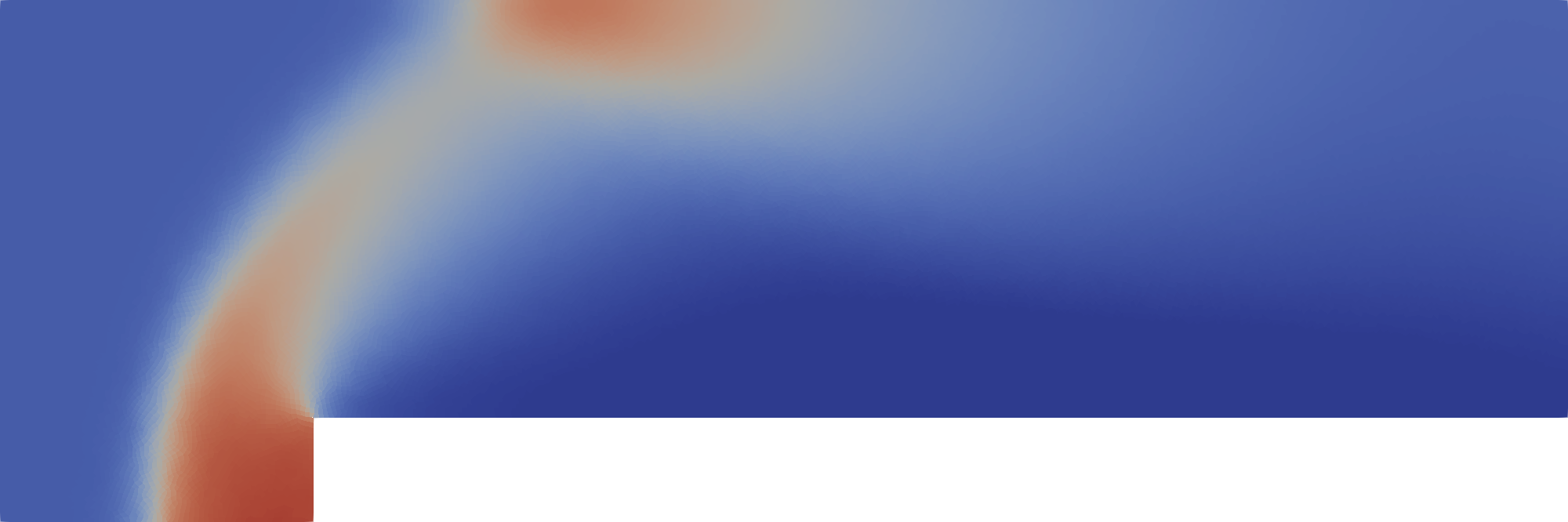}
        \put(9.4,0){\color{black}\line(0,1){33.3}}
        \put(28.7,33.3){\color{black}\line(0,-1){26.6}}
    \end{overpic}}
    \end{subfigures}
    \caption{Forward facing step horizontal velocity $\vel_x$ for $\epsilon = 10^{-2}$, $t_{\textrm{end}}=1$ for reference Euler (a), micro splitting method (b), and diffusive micro model (c). Black lines marking the bow shock and reflection point of the reference micro solution for comparison.}
    \label{fig:ffs_ref_Kn0p01}
\end{figure}

\begin{figure}[htbp!]
    \centering
    \begin{subfigures}
    \subfloat[Projective Integration, HME with $M=3$, $\delta t = 2.5\cdot 10^{-4}$, $K=1$, $\Delta t = 0.001$. \label{fig:0p01HMEM3Kn0p01deltaT0p001InnerdeltaT0p00025K1Extrapt11}
    ]{\begin{overpic}[width=0.95\textwidth]{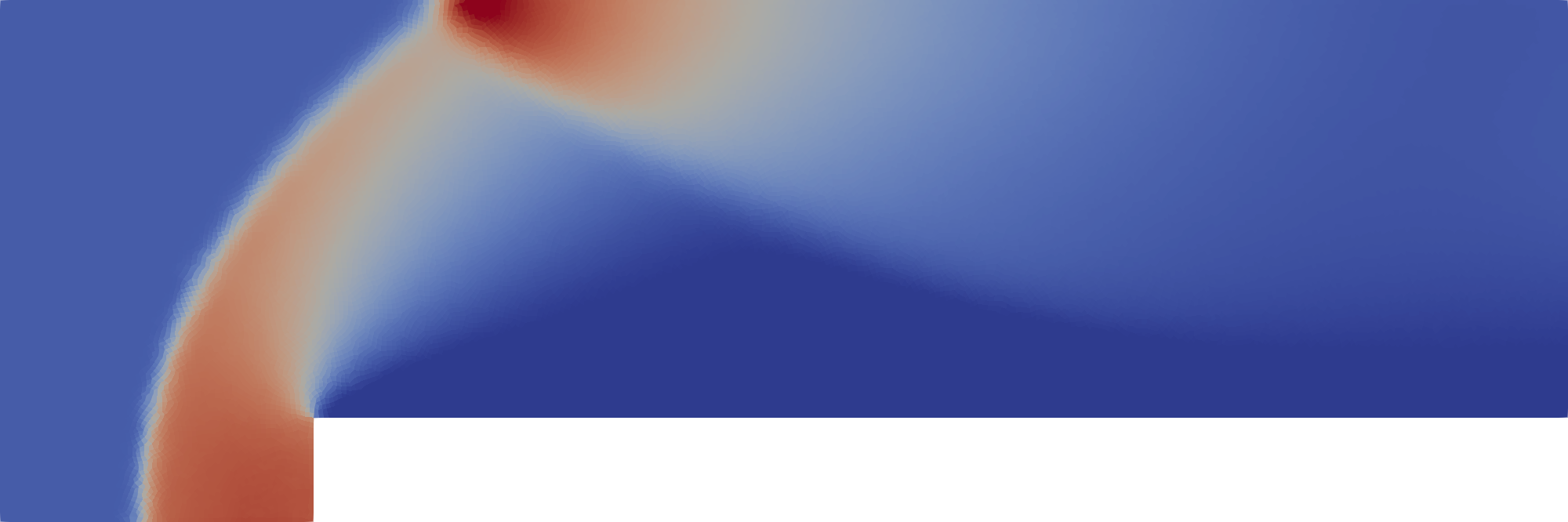}
        \put(9.4,0){\color{black}\line(0,1){33.3}}
        \put(28.7,33.3){\color{black}\line(0,-1){26.6}}
    \end{overpic}}\\
    \subfloat[Coarse Projective Integration, HME with $M=3$, $\delta t = 2.5\cdot 10^{-4}$, $K=1$, $\Delta t = 0.001$. \label{fig:0p01HMEM3Kn0p01deltaT0p001InnerdeltaT0p00025K1Coarset11}
    ]{\begin{overpic}[width=0.95\textwidth]{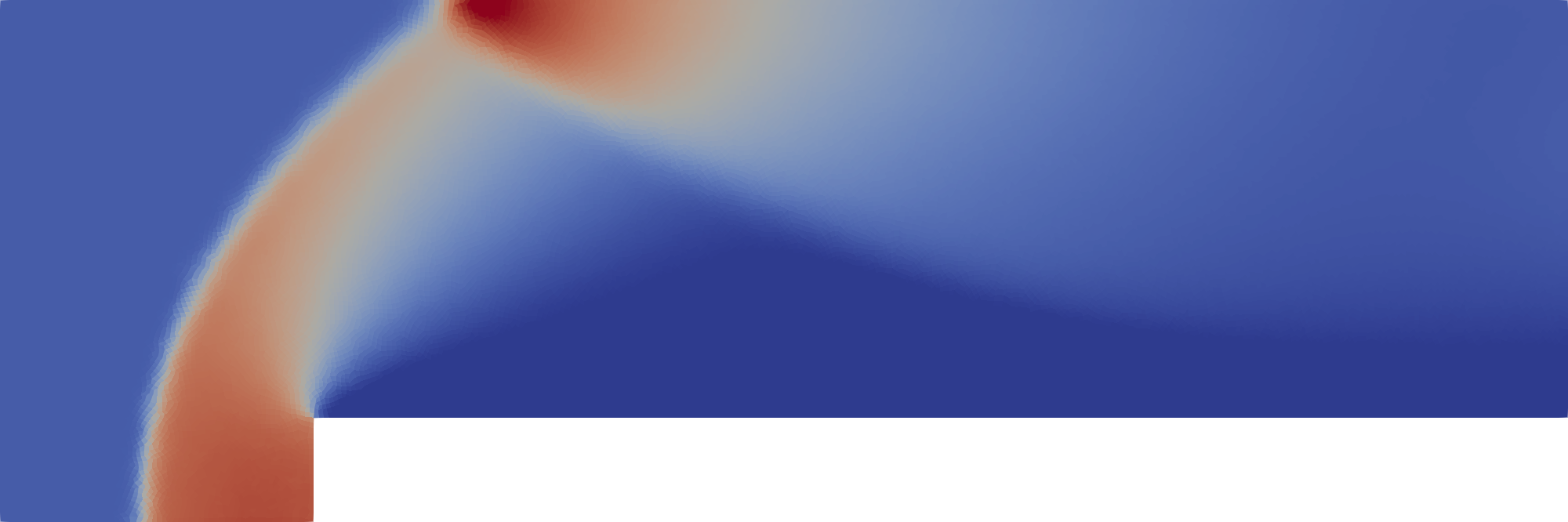}
        \put(9.4,0){\color{black}\line(0,1){33.3}}
        \put(28.7,33.3){\color{black}\line(0,-1){26.6}}
    \end{overpic}}\\
    \subfloat[Micro-macro method, HME with $M=3$, $\delta t = 2.5\cdot 10^{-4}$, $K=1$, macro Euler with $\Delta t = 0.001$. \label{fig:0p01HMEM3Kn0p01deltaT0p001InnerdeltaT0p00025K1t11}
    ]{\begin{overpic}[width=0.95\textwidth]{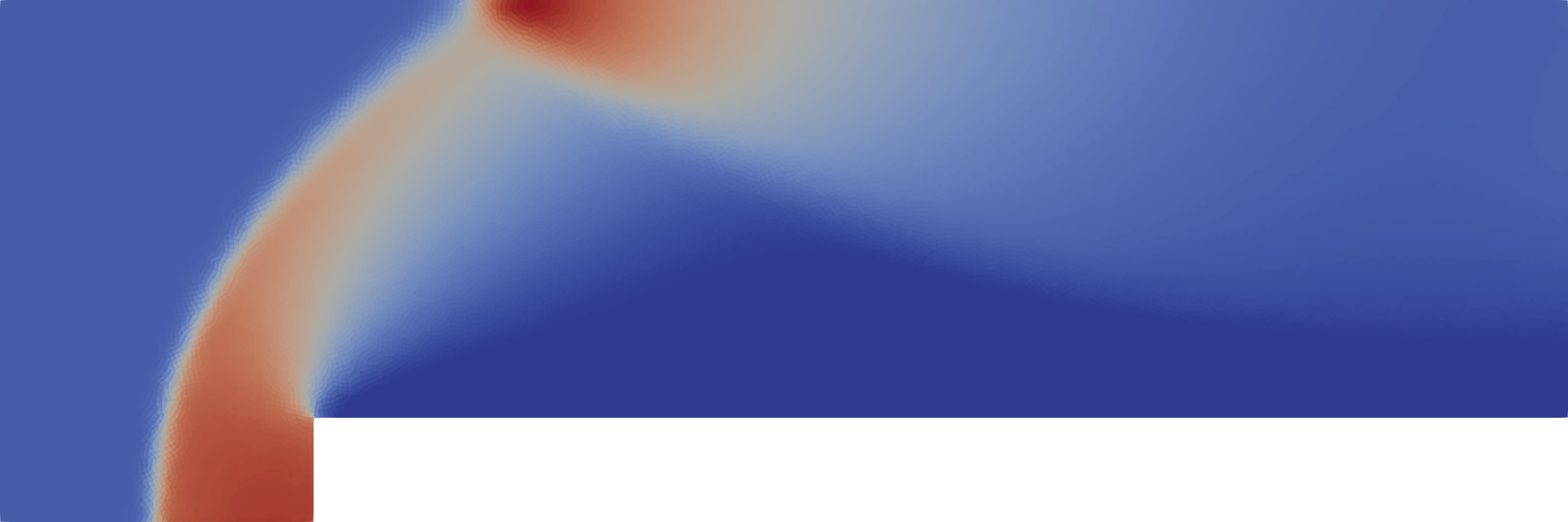}
        \put(9.4,0){\color{black}\line(0,1){33.3}}
        \put(28.7,33.3){\color{black}\line(0,-1){26.6}}
    \end{overpic}}
    \end{subfigures}
    \caption{Forward facing step horizontal velocity $\vel_x$ for $\epsilon = 10^{-2}$, $t_{\textrm{end}}=1$ for Projective Integration (a), Coarse Projective Integration (b), and micro-macro method (c). Black lines marking the bow shock and reflection point of the reference micro solution from Figure \ref{fig:ffs_ref_Kn0p01} for comparison.}
    \label{fig:ffs_Kn0p01}
\end{figure}
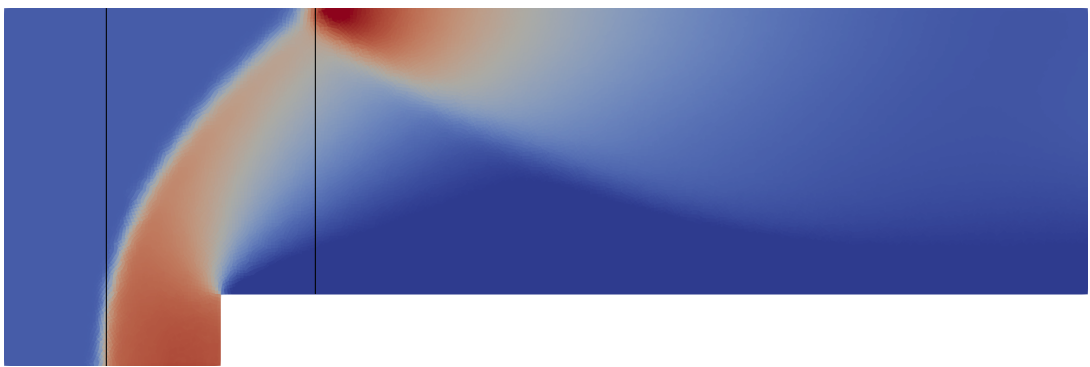
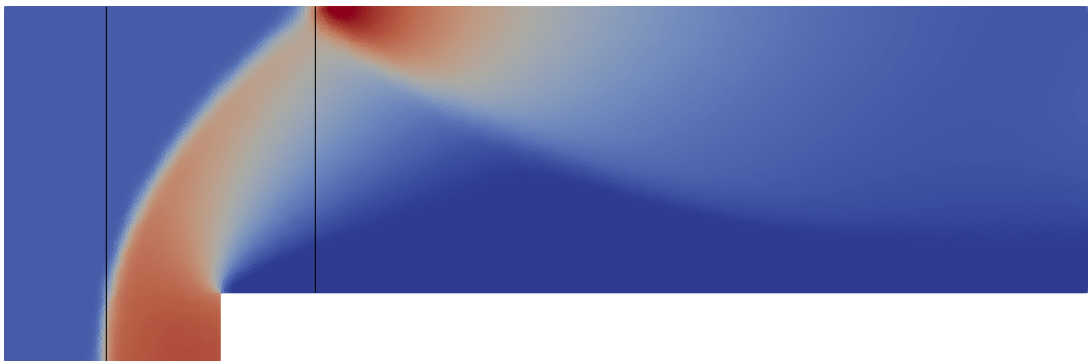
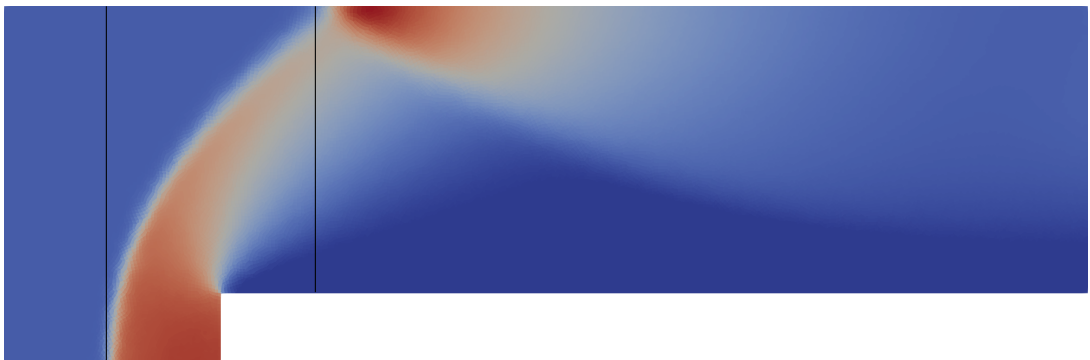

The computational grid uses $31,951$ unstructured quadrilateral grid cells of size about $\Delta x \approx 0.01$ in one direction. The 2D hierarchical micro-macro acceleration from section \ref{sec:7_mM2DHME} is used and compared with the models from \cite{Koellermeier2017d,Koellermeier2021} for reference. The macroscopic time step size according to a CFL number of $0.5$ is chosen as $\Delta t = 0.001$ and the end time is $t_{\textrm{end}}=1$ way before the steady-state so that transient phenomena can be observed and the time-accuracy can be evaluated. Figures \ref{fig:ffs_ref}-\ref{fig:ffs_Kn0p01} show the velocity in x-direction $\vel_x$, other variables show similar results.

Figure \ref{fig:ffs_ref} shows the solution of the forward facing step problem for $\epsilon = 10^{-3}$, $t_{\textrm{end}}=1$ for the Euler equations, reference micro model using simple time splitting (which is difficult to generalize or extend to higher-order, see \cite{Koellermeier2021}), and a micro model using a small time step. We see that the micro model using time splitting results in a solution very close to the Euler equations for this small relaxation time. However, the small relaxation time leads to a stiff right-hand side of the Boltzmann equation, such that an explicit micro model requires a very small time step size. In Figure \ref{fig:0p01HMEM3Kn0p001deltaT0p0001t11} this leads to excessive numerical diffusion due to the many small time steps. This adds artificial viscosity in every time step, visible by smeared out results. 

Figure \ref{fig:ffs} shows the solution of the forward facing step problem for $\epsilon = 10^{-3}$, $t_{\textrm{end}}=1$ for three of the methods outlined in Section \ref{sec:3}: The Projective Integration (PI) method from \cite{Koellermeier2021}, the Coarse Projective Integration method (CPI), and the new hierarchical micro-macro acceleration outlined in Section \ref{sec:7_mM2DHME}. Both the PI and the CPI method result in a significant offset of the shock speed, while the micro-macro method captures the position of the shock much better. Additionally, the PI method introduces excess numerical diffusion, due to the necessary $K=3$ inner iterations for this test case. The CPI as well as the micro-macro method result in a sharper shock profile. In summary, the micro-macro method yields both accurate and fast solutions to the test case.

Increasing the relaxation time $\epsilon$ towards more non-equilibrium, Figure \ref{fig:ffs_ref_Kn0p01} shows the solution of the forward facing step problem for $\epsilon = 10^{-2}$, $t_{\textrm{end}}=1$ for the Euler equations, reference micro model using time splitting, and a micro model using a small time step. The reference micro model yields a significantly different shock position and reflection position compared to the Euler solution. This is due to stronger non-equilibrium effects, that cannot be represented by the Euler equations. A standard discretization of the micro model requires a smaller time step size and results in additional numerical diffusion. This leads to a smeared out shock and hence a wrong reflection position.

Figure \ref{fig:ffs_Kn0p01} shows the solution of the forward facing step problem for $\epsilon = 10^{-2}$, $t_{\textrm{end}}=1$ for three of the methods outlined in Section \ref{sec:3}: The Projective Integration (PI) method from \cite{Koellermeier2021}, the Coarse Projective Integration method (CPI), and the new hierarchical micro-macro acceleration outlined in Section \ref{sec:7_mM2DHME}. We observe that the PI method predicts a slight offset of the shock position and the reflection position to the left whereas the CPI method seems to yield a very accurate result. The micro-macro method, on the other hand, yields a small shift of the shock position to the right and thus results in a delayed reflection point, too. Overall, the micro-macro solution has more similarity with the equilibrium Euler solution, due to the macro Euler model and the subsequent matching. Note that here the PI, CPI and micro-macro method all use only $K=1$ inner iterations, so no additional diffusion is visible. 

In summary, the hierarchical micro-macro acceleration method yields good results close to equilibrium flow. Deviations from equilibrium can be represented best with a CPI or PI method, which do not use the Euler model as macro model.

Table \ref{tab:methods_runtime_FFS} shows the measured computational speedup of the different 2D methods presented in this paper as shown in figures \ref{fig:ffs_ref} and \ref{fig:ffs} for the case $\epsilon = 10^{-3}$.
\begin{table}[h]
    \centering
    \begin{tabular}{|c||c|c|c|}
       \hline
       \bf method & \bf speedup & \bf th. speedup \\ \hline
       \hline
       Euler equations \ref{fig:0p01EulerKn0p01deltaT0p001t11scale} & 20.99 & 25 \\ \hline
       micro model \ref{fig:0p01HMEM3Kn0p001deltaT0p0001t11} & 1 & 1 \\ \hline
       PI \ref{fig:0p01HMEM3Kn0p001deltaT0p001InnerdeltaT0p0001K3Extrapt11} & 2.31 & 2.5 \\ \hline
       CPI \ref{fig:0p01HMEM3Kn0p001deltaT0p001InnerdeltaT0p0001K1Coarset11} & 4.24 & 5 \\ \hline
       micro-macro \ref{fig:0p01HMEM3Kn0p001deltaT0p001InnerdeltaT0p0001K1t11} & 3.77 & 5 \\ \hline
    \end{tabular}
    \caption{Forward facing test case speedup computed at runtime and theoretical speedup computed from section \ref{sec:Complexity} with respect to micro model for different methods from figures \ref{fig:ffs_ref} and \ref{fig:ffs}, compare Table \ref{tab:methods_runtime}.}
    \label{tab:methods_runtime_FFS}
\end{table}

It is obvious that the Euler equations from figure \ref{fig:0p01EulerKn0p01deltaT0p001t11scale} lead to the fastest solution, but they might lack accuracy in regions of existing non-equilibrium. The Projective Integration (PI) method from figure \ref{fig:0p01HMEM3Kn0p001deltaT0p001InnerdeltaT0p0001K3Extrapt11} results in a moderate speedup, which is close to the theoretical speedup expected from the complexity analysis in section \ref{sec:Complexity}. The Coarse PI (CPI) method results in a larger speedup when compared to the hierarchical micro-macro acceleration, which is due to the simpler extrapolation instead of a physical Euler model. The use of a physical macro model reduces the speedup of the hierarchical micro-macro acceleration, but might increase its accuracy beyond a simple extrapolation of the solution. Note that the speedup of both the CPI scheme and the hierarchical micro-macro acceleration differs by some fraction from the theoretical speedup computed in section \ref{sec:Complexity} because the matching, restriction, and macro steps were neglected for the theoretical computation. However, the new hierarchical micro-macro acceleration still yields a significant speedup for this test case.
\section{Conclusion}
\label{sec:conclusion}

In this paper, we introduced a hierarchical micro-macro acceleration for moment models of kinetic equations. Variants of the hierarchical micro-macro acceleration promise fast computation with similar accuracy compared to a full micro model or better accuracy than a macro model with still reasonable computation time.
The hierarchical micro-macro acceleration is based subsequently on a micro step, a restriction to the macroscopic variables, a macro step, and a matching step to reconstruct the microscopic variables. There is a large amount of flexibility to construct these four steps, so that many existing methods like Projective Integration (PI) and Coarse Projective Integration (CPI) and many new micro-macro schemes can be easily compared to the hierarchical micro-macro acceleration. While a very general matching operator can be used, we employ $L^2$ matching, which allows for a fast solution of the matching problem. 
Numerical results of a 1D shock structure test case and a 2D forward facing step application case showed good accuracy with a speedup of up to 250 for our new micro-macro method.

The research in this paper opens up possibilities for many extensions. Different moment models like maximum entropy models or regularized moment models can be readily used for the micro/macro steps. Different matching operators like relative entropy can improve the physical accuracy of the matching step. Finally, an investigation of convergence properties and error estimates for the whole method would be desirable.

\section*{Data Availability Statement}
The datasets generated and analysed during this study are available
from the corresponding author on reasonable request.

\section*{Acknowledgements}
This research has been partially supported by the European Union’s Horizon 2020 research and innovation program under the Marie Sklodowska-Curie grant agreement no. 888596. 
The authors would like to acknowledge the financial support of the CogniGron
research center and the Ubbo Emmius Funds (University of Groningen). 




\bibliographystyle{abbrv}
\bibliography{paper}

\end{document}